\documentclass[10pt]{amsart}
\allowdisplaybreaks[4]
\usepackage{amssymb,color}
\usepackage[square,compress,comma, numbers,sort]{natbib}
\usepackage[colorlinks=true, citecolor=blue, linkcolor=blue]{hyperref}
\usepackage{amsthm}
\usepackage{mathtools}

\newcommand{\E}[1]{\mathbb{E}\left\{ #1\right\}}

\newcommand{\pk}[1]{\mathbb{P} \left(#1 \right) }

\definecolor{c20}{rgb}{0.,0.7,0.}
\definecolor{c30}{rgb}{0.,0.,1.}
\definecolor{c40}{rgb}{1,0.1,0.7}
\definecolor{c50}{rgb}{1,0,0}
\definecolor{c60}{rgb}{0,0.9,0.1}

\def\td{\text{\rm d}}

\newcommand{\abs}[1]{\lvert #1 \rvert}
\newcommand{\ABs}[1]{ \biggl \lvert #1 \biggr \rvert}
\newcommand{\norm}[1]{\lVert #1 \rVert}

\newcommand{\R}{\mathbb{R}}

\newcommand{\inr}{\in \R}

\newcommand{\limit}[1]{\lim_{#1 \to   \infty}}

\newcommand{\inv}[1]{\overleftarrow{#1}}

\newcommand{\BQN}{\begin{eqnarray}}
\newcommand{\EQN}{\end{eqnarray}}
\newcommand{\BQNY}{\begin{eqnarray*}}
\newcommand{\EQNY}{\end{eqnarray*}}
\def\nncol#1{\textcolor{black}{#1}}

\def\K1#1{\textcolor{cyan}{#1}}
\def\K1#1{#1}

\def\kk#1{\textcolor{cyan}{#1}}

\def\bqny#1{ \nncol{ \begin{eqnarray*} #1 \end{eqnarray*}}}
\def\bqn#1{ \nncol{ \begin{eqnarray} #1 \end{eqnarray}}}

\newcommand{\BS}{\begin{sat}}
\newcommand{\ES}{\end{sat}}
\newcommand{\BT}{\begin{theo}}
\newcommand{\ET}{\end{theo}}
\newcommand{\BK}{\begin{korr}}
\newcommand{\EK}{\end{korr}}

\newcommand{\BD}{\begin{de}}
\newcommand{\ED}{\end{de}}

\newcommand{\BIT}{\begin{itemize}}
\newcommand{\EIT}{\end{itemize}}
\newcommand{\BDI}{\begin{description}}
\newcommand{\EDI}{\end{description}}

\newcommand{\BRM}{\begin{remark}}
\newcommand{\ERM}{\end{remark}}

\newcommand{\BEL}{\begin{lem}}
\newcommand{\EEL}{\end{lem}}

\newtheorem{theo}{Theorem}[section]
\newtheorem{sat}[theo]{Proposition}
\newtheorem{de}[theo]{Definition}
\newtheorem{lem}[theo]{Lemma}

\newtheorem{korr}[theo]{Corollary}
\newtheorem{remark}[theo]{Remark}

\newcommand{\netheo}[1]{{Theorem \ref{#1}}}

\newcommand{\prooftheo}[1]{ \textbf{Proof of Theorem} \ref{#1} }

\newcommand{\prooflem}[1]{\textbf{Proof of Lemma} \ref{#1}}
\newcommand{\proofkorr}[1]{\textbf{Proof of Corollary} \ref{#1}}

\newcommand{\COM}[1]{}

\newcommand{\QED}{\hfill $\Box$ \\}

\def\rw{\rightarrow}

\def\IF{\infty}

\def\LT{\left}
\def\RT{\right}

\def\eHH#1{\textcolor{c50}{#1}}

\def\eHH#1{#1}

\def\ehd#1{\textcolor{c50}{#1}}
\def\ehd#1{#1}
\def\ehe#1{\textcolor{black}{#1}}

\def\blu#1{\textcolor{black}{#1}}
\def\K#1{\textcolor{black}{#1}}
\def\kk#1{\textcolor{black}{#1}}
\def\rrd#1{\textcolor{black}{#1}}
\def\re#1{\textcolor{black}{#1}}

\def\vf{\sigma^2}

\def\B{\mathcal{B}}
\def\k2#1{\textcolor{cyan}{#1}}
\begin{document}

\title
{  Sojourn times of Gaussian processes with trend}

\author{Krzysztof D\c{e}bicki}
\address{Krzysztof D\c{e}bicki, Mathematical Institute, University of Wroc\l aw, pl. Grunwaldzki 2/4, 50-384 Wroc\l aw, Poland}
\email{Krzysztof.Debicki@math.uni.wroc.pl}
\author{Peng Liu}
\address{Peng Liu, Department of Actuarial Science, Faculty of Business and Economics, University of Lausanne, UNIL-Dorigny 1015 Lausanne, Switzerland and Department of Statistics and Actuarial Science, University of Waterloo, Canada}
\email{peng.liu1@uwaterloo.ca}
\author{Zbigniew Michna}
	\address{Zbigniew Michna, Department of Mathematics and Cybernetics, Wroc\l aw University of Economics, Poland}
	\email{zbigniew.michna@ue.wroc.pl}

\bigskip

\date{\today}
 \maketitle
\bigskip

\begin{abstract}
We derive exact tail asymptotics of sojourn time
above
the level \re{$u\geq 0$}
\[
\pk{v(u)\int_0^T \mathbb{I}(X(t)-ct>u)\td t>x}, \quad x\geq 0
\]
as $u\rw\IF$, where $X$ is a Gaussian process with continuous sample paths,
$c$ is some constant,  $v(u)$ is a positive function of $u$ and $T\in (0,\IF]$.
Additionally, we analyze asymptotic distributional properties of
\[\tau_u(x):=\inf\left\{t\geq 0: \ehe{v(u)}
\int_0^t \mathbb{I}(X(s)-cs>u)\td s>x\right\},
\]
as $u\to\infty$, $x\geq 0$, where $\inf\emptyset=\IF$.
The findings of this
contribution are illustrated by a detailed analysis
of a class of Gaussian processes with stationary increments
{and a family of self-similar Gaussian processes.}
\end{abstract}



\section{Introduction}
Let $Y(t) ,t\inr$ be a centered random process with \re{c\`{a}dl\`{a}g} sample paths and let for $T>0$
\K{
\begin{eqnarray}
\int_0^T \mathbb{I}_u(Y(t)) \td t \label{m1}
\end{eqnarray}
be the sojourn time of $Y$  above the level $u\inr$ in interval $[0,T]$,
where
$\mathbb{I}_u(x)=\mathbb{I}(x>u)$.
}
The \K{asymptotic properties of
(\ref{m1}), as $u\to\infty$ for $Y$ being a centered Gaussian process}
has been extensively studied by Berman, see e.g., \cite{Berman82,MR803245,berman1987extreme, Berman92}.
\re{An} explicit form of the distribution of (\ref{m1}) is known only for very few special processes.
In particular, for
$Y=B_1$ being a standard Brownian motion,
by the arcsin law of Paul  L\'{e}vy, we have
\BQN
\pk{\int_0^T \mathbb{I}_{\ehd{0}}(B_1(t)) \td t  > x} =
1- \frac{2}{\pi}\arcsin\Bigl( \sqrt{\frac{x}{T}}\Bigr), \quad 0 < x < T < \IF.
\EQN
An extension of this arcsin law \K{is} obtained in \cite{AKA} for the case \K{of}
$Y(t)= B_1(t) - ct, c\not=0$. For the infinite time horizon, i.e., $T=\IF$
and $c>0$,
in view of  \cite{Borodin2002}[Eq. (3), p. 255]
$$\pk{\int_0^\infty \mathbb{I}_0(B_1(t)-ct) \td t  \in \td y}=\left(\frac{\sqrt{2}c}{\sqrt{\pi y}}e^{-\frac{c^2y}{2}}-\frac{2c^2}{\sqrt{\pi}}\int_{\frac{c\sqrt{y}}{\sqrt{2}}}^\IF e^{-v^2}\td v\right)\td y, \quad y\geq 0,$$
which implies that  for any $x\geq 0$
\BQN\label{bermanbrownain}
\pk{\int_0^\infty \mathbb{I}_0(B_1(t)-ct) \td t>x}=2(1+c^2x)\Psi(c\sqrt{x})-\frac{c\sqrt{2x}}{\sqrt{\pi}}e^{-\frac{c^2x}{2}},
\EQN
where $\Psi$ denotes the \K{tail} distribution of a standard normal random variable $N(0,1)$.

\K{The above result can be extended to the \re{class of} spectrally negative
L\' evy \re{processes} with a negative drift and a general level $u$.
Indeed, let $X$ be a spectrally negative L\'{e}vy process.
Then for any non-negative $\lambda$
 $$\E{ e^{\lambda X(t)}}=e^{t \psi(\lambda)}, \quad t\ge 0, $$
with $\psi$ being a strictly convex function such that
$\limit{\lambda} \psi(\lambda)= \IF$ and $\psi'(0+)=\E{X(1)}$.
\BT \label{ThZ} Suppose that $X(t),t\ge 0  $ is a spectrally negative L\'{e}vy process
  such that   $\E{X(1)}<c$ for some $c>0$. Then \ehe{for any \re{$u\geq 0$}}
 \BQN \label{eq:specn}
 \pk{\int_0^\IF \mathbb{I}_u(X(t)-ct)\td t>x}=e^{-\alpha u}\,\pk{\int_0^\IF \mathbb{I}_0(X(t)-ct)\td t>x},
 \EQN
 where $\alpha>0$ is the unique  positive solution \re{to} $ \psi(\alpha)=c\alpha$.
 \ET
\re{In} the particular case that $X=B_1$, we have $\E{B_1(1)}=0< c$ and further
$\psi(\alpha)=\frac{\alpha^2}{2}=c\alpha$
has the unique positive solution $\alpha=2c$. Consequently, for any $x\geq 0$ and \ehe{$u\inr$}
\netheo{ThZ} combined with (\ref{bermanbrownain}) implies
\BQN\label{brownian}
 \pk{\int_0^\infty \mathbb{I}_u(B_1(t)-ct)\td t>x}=\left(2(1+c^2x)\Psi(c\sqrt{x})-\frac{c\sqrt{2x}}{\sqrt{\pi}}e^{-\frac{c^2x}{2}}\right)e^{-2c u}.
 \EQN
}

The study of distributional properties of
\K{occupation-type functionals for L\'evy processes}
is  crucial for many applications in finance and insurance (e.g. the occupation time in red or the inverse occupation time - the time of cumulative Parisian ruin), see for instance \cite{AKA,OKP3,OKP4, OKP7,AMINE19}. The number of papers dealing with occupation times (sojourn times) is huge; most of the articles discuss the derivation of
Laplace transform, see the recent contributions \cite{OKP1,OKP2, OKP5,OKP6}
and \kk{references} therein. {Recent paper \cite{AMINE19} derives the density of occupation time for spectrally
negative L\'{e}vy processes with exponential time horizon.}

In this paper we shall focus on \K{analogues of Theorem \ref{ThZ} for a wide
{class of Gaussian processes. }}
Since
\ehe{the distribution \kk{of} (\ref{m1})
is not tractable in the general Gaussian case, our investigation is
concerned with the derivation of the \kk{exact} asymptotic behavior of}
\BQN\label{eq0}
\ehd{p_{T}(u,x):}=\pk{v(u)\int_0^T \mathbb{I}_u(X(t)-ct)\td t>x}, \quad x\geq 0
\EQN
as $u\rw\IF$, where  $X$ is a Gaussian process with continuous sample paths,
$v(u),u>0$ is a \ehd{positive}
function (will be specified below) and $T\in (0,\IF]$.
{In order to avoid trivialities, we suppose that $c>0$ if $T=\infty$, while for $T<\infty$ we allow $c\in\R$.}
In contrast to the classical results
for centered Gaussian processes $X$ by Berman, see e.g.,
\cite{Berman82,MR803245,berman1987extreme, Berman92},
where the asymptotics
of
$\pk{v(u)\int_0^T \mathbb{I}_u(X(t))\td t>x}$ as $u\to\infty$
 \re{is} given for a.e. $x\geq 0$,
an important advantage of the technique used in this contribution
allowed us to show that the  asymptotics of (\ref{eq0}) \re{holds} for all $x\geq 0$.
This is due to continuity of Berman-type constant which shall be proven in Lemma \ref{Pickands}.

\kk{Additionally, motivated by recent \kk{investigations on the ruin time, see e.g., \cite{HP2008} and \cite{HJ2014}, }
we shall analyze asymptotic distributional properties of}
\BQN\label{passengetime}\tau_u(x):=\inf\left\{t\geq 0: \ehe{v(u)}
\int_0^t \mathbb{I}_u(X(s)-cs)\td s>x\right\},
\EQN
as $u\to\infty$, $x\geq 0$,
with the convention that $\inf\emptyset=\IF$.
\rrd{Note that  $\tau_u(x)$ is called  the cumulative Parisian ruin time in \cite{OKP7}, which is of interest in risk theory;
 $\tau_u(0)$ is the first passage time \re{of the level $u$ by the process $X(t)-ct$},
  which is referred to as {{\it the ruin time}} in \cite{HP2008}.}
 In Section 3 a distributional approximation, as $u\to \IF$, of  $\tau_u^*(x_1, x_2)$ defined by
\BQN\label{conditionaltime}\tau_u^*(x_1,x_2):=\tau_u(x_2)\Big {|}\tau_u(x_1)<\IF, \quad 0\leq x_1\leq x_2<\IF
\EQN
is derived.

Brief outline of the paper:
\K{
In Section \ref{s.prel} we derive \re{an} exact approximation of sojourn times for general Gaussian processes. In Section \ref{s.examples}
we apply this result  to
(\ref{eq0}),  for $X$ being a centered  Gaussian process with stationary increments and a self-similar Gaussian process.
Both scenarios $T\in (0,\infty)$ and $T=\infty$ are considered.
Section \ref{s.lemmas} contains some lemmas that are useful in the proofs of the main results,
while the proofs are presented in Section \ref{s.proofs}.
}

\section{Main result}\label{s.prel}
{In this section we provide a general result preparing us to derive exact
asymptotics of (\ref{eq0})
for a wide class of
centered Gaussian process $X(t),t\inr$ with continuous trajectories.
}
In order to motivate the study of this section write first for $u>0,x\ge 0, {c>0}$  and $v(u)$  representing an arbitrary positive function
\BQN\label{eq1}
p_\IF(u,x)&:=&\pk{v(u)\int_0^{\infty} \mathbb{I}_u(X(t)-ct)\td t>x}\notag\\
&=&\pk{v(u)u\int_0^{\infty} \mathbb{I}_0(X(ut)-cut-u)\td t>x}\nonumber\\
&=&
\K{\pk{v(u)u\int_0^{\infty} \mathbb{I}_{M(u)}\left(\frac{X(ut)}{u(1+ct)}M(u)\right)\td t>x}},
\EQN
where
$M(u)=\inf_{t\in \K{[0,\IF)}}\frac{u(1+ct)}{\sqrt{Var(X(ut))}}$
\ehd{\kk{is} assumed to be positive.} {Then} $\frac{X(ut)}{u(1+ct)}M(u)$ is a Gaussian process with mean $0$ and $\sup_{t\geq 0}Var\left(\frac{X(ut)}{u(1+ct)}M(u)\right)$=1. \K{As it will be proven in Section \ref{s.proofs}, for
	$\delta_u$, a properly} chosen function of $u$ with $\lim_{u\rw\IF}\delta_u=0$, and  $t_u=\arg \inf_{t\geq 0}\frac{u(1+ct)}{\sqrt{Var(X(ut))}} $
for $x\ge0$ we have
\BQNY
p_\IF(u,x)& \sim &  \pk{v(u)u\int_{[t_u-\delta_u, t_u+\delta_u]} \mathbb{I}_{\K{M(u)}}\left(\frac{X(ut)}{u(1+ct)}M(u)\right)\td t>x}\\
&=&\pk{\int_{[-v(u)u\delta_u, v(u)u\delta_u]} \mathbb{I}_{\K{M(u)}}\left(\frac{X(ut_u+t/v(u))}{u(1+ct_u)+\K{c}t/v(u)}M(u)\right)\td t>x}.
\EQNY

\ehd{A similar transformation for $p_T(u,x), T \in (0,\IF)$ 
shows that \kk{in general the problem to deal with can be reduced to}
 }
\BQN\label{mainasym0}
p_T(u, x)\sim \pk{\int_{E(u)} \mathbb{I}_{n(u)}(Z_u(t) ) \td t>x}, \quad u\rw\IF, \EQN
where $x\ge 0$, $n(u)$ is a function of $u$ and  $Z_u,u>0$ is a family of centered Gaussian processes with continuous trajectories  defined  on the 
interval $E(u)=[a_1(u), a_2(u)]$.


\COM{
Note that for $x=0$, (\ref{eq1}) has been extensively studied in the literature for $X$ being a self-similar process and (OR?) a Gaussian process with stationary increments, see e.g., \cite{HP99}, \cite{DE2002}, \cite{HP2004} and \cite{DI2005}. In the aforementioned papers,  we observe that
\BQNY
\pk{\sup_{t\geq 0} \frac{X(ut)}{u(1+ct)}n(u)>n(u)}\sim
 \pk{\sup_{t\in (t_u-\delta_u, t_u+\delta_u)} \frac{X(ut)}{u(1+ct)}n(u)>n(u)},
\EQNY
}
\kk{In the rest of this paper we
shall impose some standard} assumptions on the behaviour of the variance function $\sigma_u$ and \K{the} correlation function $r_u$ of $Z_u$.
\kk{In particular}, we shall assume that 
\BQN\label{var}
\ehd{\limit{u} \sup_{t\neq 0, t\in E(u)}\ABs{ \frac{ 1-\sigma_u(t)}{ w(g(u)|t|)}-1}=0},
\EQN
where $w$ is a \K{positive} regularly varying function at $0$ with index $\beta>0$, and $g(u)$ \K{satisfies} $\lim_{u\rw\IF}g(u)=0$.

\ehd{In the following
	$\Delta(u),n(u), u>0$ are positive \blu{functions} such that
	\begin{eqnarray}\label{delta(u)}
	\lim_{u\rw\IF}\Delta(u)=\varphi\in [0,\IF], \quad
	\limit{u} n(u)= \IF.
	\label{phi1}
	\end{eqnarray}
For the correlation function of $Z_u$ we shall assume that for  $\Delta(u)$ satisfying \eqref{phi1} we have}
\BQN\label{cor}
\lim_{u\rw\IF}\sup_{s\neq t, s,t\in E(u)}\left|\frac{n^2(u)\left(1-r_u(s,t) \right)}{\frac{\sigma^2_\eta(\Delta(u)|t-s|)}{\sigma^2_\eta(\Delta(u))}}-1\right|=0,
\EQN
where
$\eta(t),t\inr$ is a centered  Gaussian process with continuous trajectories,
stationary increments and variance function
$\sigma^2_\eta(t)>0, t>0,$ being
regularly varying at $0$ and \kk{at} $\IF$ with indexes
\kk{$2\alpha_0\in (0,2]$ and $2\alpha_\IF\in (0,2)$}, respectively.

\K{Assumptions (\ref{var}) and (\ref{cor})
are satisfied for large classes of Gaussian processes, see e.g.,
\cite{HP99}, \cite{DE2002}, \cite{HP2004} and \cite{DI2005}.}
\kk{For example, they}  \re{are compatible with} those in Theorem 3.2 in \cite{KEP2016}.\\		
\ehd{Next, for any $\varphi\in [0,\IF]$ set}
\BQNY
\eta_{\varphi}(t)=\left\{\begin{array}{cc}
	B_{\re{2}\alpha_0}(t),& \varphi=0\\
	\frac{\eta(\varphi t)}{\sigma_\eta(\varphi)}, & \varphi\in (0,\IF)\\
	B_{\re{2}\alpha_\IF}(t), & \varphi=\IF\,,
\end{array}
\right.
\EQNY
\kk{where $B_\alpha$ is a fractional Brownian motion (fBm)
with self-similarity index $\alpha/2 \in (0,1]$
and process $\eta$ is defined above.}
\\
\ehd{For a random process $W(t),t\inr$ with continuous trajectories},  $x\geq0$,  $E$ a compact subset of $\mathbb{R}$ and $h$ a continuous function on $E$  we define
\BQNY
\B_{W}^{h}(x,E)&=&\int_{\mathbb{R}}\pk{\int_{E} \mathbb{I}_0(\sqrt{2}W(t)-Var(W(t))-h(t)+z)\td t>x}e^{-z}\td z
\EQNY
and when the limit exists, we set
$$\ehd{\B_{W}^h(x)=\lim_{S\rw\IF}\frac{\re{\B_{W}^h}(x,[0,S])}{S^{\mathbb{I}(h=0)}}.}$$
Moreover, set below
 $$\widehat{\B}_{B_\alpha}^{h}(x)=\lim_{S\rw\IF}\B_{B_\alpha}^{h}(x, [-S,S])\,,$$
 provided that the above  limit is finite.
\kk{For $h= 0$, we suppress the superscript and write  $\B_{W}(x)$ or $\B_W(x,E).$}
If $W=B_\alpha$ with $B_\alpha$ fractional Brownian motion,
 then \re{$\B_{B_\alpha}(0)$} is simply the Pickands constant, see e.g.,\cite{Pit96}, \cite{DE2002}, \cite{HP2004} and \cite{DI2005}.\\
  If $h$ is strictly positive, then
$  \B_{B_\alpha}^{h}(0)$  and $\widehat{\B}_{B_\alpha}^h(0)$ reduce to Piterbarg constants.
{We refer to} \cite{HP99,DE2002, HP2004, DI2005} for the existence and properties of  Pickands and Piterbarg related constants. \ehe{All our asymptotic results below  hold for all $x\geq0$, which generalize the results in \cite{Berman82}-\cite{Berman92} and \cite{KEZX17}, where the asymptotics hold for almost all $x\in [0,\IF)$.}

\ehe{Throughout this paper,  $\inv{f}$ stands for the generalized asymptotic (unique)
inverse of a regularly varying function $f$; see \cite{BI1989}}.

Next we present \kk{the main result of this contribution}.
\def\THU{\ehe{\theta(u)}}
\BT\label{TH1}
{Let $Z_u(t),  t\in E(u)$ with $ E(u):=[a_1(u), a_2(u)]$ be a family of centered Gaussian processes with continuous
trajectories. Suppose that (\ref{var})-(\ref{cor}) hold, $\THU=\inv{w}(n^{-2}(u)) / g(u) $} and
\BQN\label{gammalim}\lim_{u\rw\IF}n^2(u)w(g(u))=\gamma \in [0,\IF], \quad \lim_{u\rw\IF}g(u)|a_i(u)|=0, i=1,2.
\EQN
i) If  $\gamma=0$ and
$$\lim_{u\rw\IF}n^2(u)w(g(u)|a_i(u)|)=x_i,
 \quad \lim_{u\rw\IF}n(u)w(g(u)|a_i(u)|)=0,  \quad i=1,2,$$  then
\BQN \label{madr1}
\pk{\int_{E(u)} \mathbb{I}_{n(u)}(Z_u(t) ) \td t>x}
\sim\B_{\eta_{\varphi}}(x)\frac{1}{\beta}\int_{y_1}^{y_2}|t|^{1/\beta-1}e^{-|t|}\td t\,
\THU 
\Psi(n(u)),
\EQN
with   $y_2-y_1>0$ and
\BQN\label{y112}
y_i=x_i \mathbb{I} ( x_i>0, \lim_{u\rw\IF} a_i(u)=\IF\})-x_i \mathbb{I}( \re{x_i>0}, \lim_{u\rw\IF}a_i(u)=-\IF), \quad i=1,2. \EQN
ii) If  $\gamma\in (0,\IF)$ and $\lim_{u\rw\IF}a_i(u)=a_i, i=1,2,$  with $a_1\in [-\IF, 0]$, $a_2\in [0,\IF]$ and $a_2-a_1>x$, then
\BQN \label{madr2}
\pk{\int_{E(u)} \mathbb{I}_{n(u)}(Z_u(t) ) \td t>x}
 \sim \B_{\eta_{\varphi}}^{\gamma|t|^\beta}(x,[a_1,a_2])\Psi(n(u)).
\EQN
iii) If  $\gamma=\IF$ and
$$\lim_{u\rw\IF}\frac{a_i(u)}{\THU} 
=b_i, \quad  i=1,2,
$$  with $b_1\in [-\IF, \IF), b_2\in (-\IF, \IF],
$ and $b_2-b_1>x,$   then 
\BQN \label{madr3}
\pk{\int_{E(u)} \mathbb{I}_{n(u)}(Z_u(t) ) \td t>\THU x}
\sim \ehe{\mathcal{B}_0}^{|t|^\beta}(x,[b_1,b_2])\Psi(n(u)).
\EQN
\ET
\BRM\label{remark} i) If we assume that $\lim_{u\rw\IF}\Delta(u)=0$ and $\sigma^2_\eta$ in (\ref{cor}) is a non-negative  regularly varying function at $0$ with index $2\alpha_0\in (0,2]$,  then Theorem \ref{TH1} still hold with $\eta_{\varphi}$ replaced by $B_{2\alpha_0}$.\\
ii) The {case $x=0$ in  Theorem \ref{TH1} generalizes}
the results in \cite{HP99}, \cite{DE2002}, \cite{HP2004} and \cite{DI2005} and together with i)
of this remark covers the results for one dimensional case in  \cite{Pit96}.
\ERM
In the following lemma we calculate the exact value of {two} special Berman constants.
\BEL\label{Constants}
 For $\gamma, \beta>0$
\BQN\label{Piterbarg}\mathcal{B}_0^{\gamma t^\beta}(x)=\mathcal{B}_0^{\gamma t^\beta}(x, [0,\IF))=e^{-\gamma x^\beta},
\quad 
\ehe{\mathcal{B}_0}^{\gamma |t|^\beta}(x, (-\IF,y])=\left\{\begin{array}{cc}
e^{-\gamma(x-y)^\beta}, & y<x/2\\
e^{-\gamma2^{-\beta} x^\beta}& y\geq x/2,
\end{array}\right.
\EQN
\ehe{and  $\ehe{\mathcal{B}_0}^{\gamma|t|^\beta}(x, [b_1,b_2])$ is continuous for all $x\ge 0$.}
{For any $x\geq 0$ and $\gamma>0$
\BQN\label{Pit2}
\widehat{\B}_{B_{2}}^{\gamma t^2}(x)=\sqrt{\frac{1+\gamma}{\gamma}}e^{-\frac{(1+\gamma)x^2}{4}}.
\EQN}
\EEL

\section{Applications}\label{s.examples}
In this section we apply Theorem \ref{TH1} to two classes of Gaussian processes:
Gaussian processes with stationary increments and self-similar Gaussian processes.
\subsection{Gaussian processes with stationary increments}\label{s.sta}
\ehe{Given $X(t),t\inr$ a Gaussian process with stationary increments and continuous sample paths, we} \K{ consider}
\BQN\label{PT}
\ehe{p_T(u,x)}\coloneqq \pk{v(u)\int_0^T \mathbb{I}_u(X(t)-ct)\td t>x}, \quad x\geq 0,
\EQN
where  the positive scaling function $v$ will be specified later, and $T\in (0,\IF]$.
Let  $\sigma^2(t)=Var(X(t))$. \K{We distinguish two cases, leading to qualitatively different asymptotics, namely  $T\in (0,\IF)$
and $T=\IF$.}
\subsubsection{Infinite time horizon.}\label{s.sti}
\K{By the stationarity of increments of $X$} the covariance function of $X$ is completely determined by its variance function $\sigma^2$.
Along the same lines as in \cite{DI2005} or \cite{KrzysPeng2015}, we assume that \\
\\
{\bf AI}: $\vf(0)=0$ and $\vf(t)$ is regularly varying at $\IF$ with index $2\alpha_\IF\in(0,\eHH{2)}$.
Further, $\vf(t)$ is twice continuously differentiable on $(0,\IF)$ with its first derivative
$\dot{\vf}(t)\coloneqq\frac{{\rm d} \sigma^2}{{\rm d}t}\left(t\right)$
and second derivative
$\ddot{\vf}(t)\coloneqq\frac{{\rm d^2} \sigma^2}{{\rm d}t^2}\left(t\right)$

being ultimately monotone at $\IF$.\\
{\bf AII}: $\vf(t)$ is regularly varying at $0$ with index $2\alpha_0\in(0,2]$.\\
\\
\K{Assumptions {\bf AI-AII} cover  a wide range of Gaussian processes with stationary increments, including two important families:
(1) {\it fractional Brownian motions} $B_{\alpha}(t)$, $\alpha\in (0,2]$
and (2) {\it Gaussian integrated processes}, i.e., the case where
$X(t)=\int_0^tZ(s)\td s,$ with $Z$ a centered continuous stationary Gaussian process with variance $1$ and correlation function {$r(s)=Cov(X(t), X(t+s)), s, t\geq 0$} satisfying some regularity conditions; see, e.g., \cite{DE2002}, \cite{HP2004}, \cite{DI2005}, \cite{KrzysPeng2015} and \cite{KEP2015}.
}
\\
{Suppose that $c>0$ and}
\ehd{let in the following ($\overleftarrow{\sigma}$ stands for the asymptotic inverse of $\sigma$)}
\BQN \label{zhr}
1/v(u)=
\overleftarrow{\sigma}
 \left(\frac{\sqrt{2}\sigma^2(ut^*)}{u(1+ct^*)}\right), \quad t^*=\frac{\alpha_\IF}{c(1-\alpha_\IF)}.
 \EQN
According  to (\ref{eq1}), recall that
\BQN\label{scale}
p_\IF(u,x)
=\pk{u v(u)\int_0^\IF \mathbb{I}_{\K{M(u)}}\left(\frac{X(ut)}{u(1+ct)}M(u)\right)\td t>x},
\EQN
where $$M(u)=\inf_{t>0}Var^{-1/2}\left(\frac{X(ut)}{u(1+ct)}\right)=\inf_{t>0}\frac{u(1+ct)}{\sigma(ut)}.$$
\K{We note} that
$$\frac{X(ut)}{u(1+ct)}M(u), \quad t\geq 0$$
 is a family of  centered Gaussian processes with \re{the} maximum of \re{their} variance \re{functions} equal to $1$. \\
Applying Theorem \ref{TH1} we arrive at the following results, where
\BQN\label{AB}
A=\left(\frac{{\alpha_\IF}}{c(1-{\alpha_\IF})}\right)^{-\alpha_\IF}\frac{1}{1-\alpha_\IF}, \quad
B=\left(\frac{\alpha_\IF}{c(1-\alpha_\IF)}\right)^{-\alpha_\IF-2}\alpha_\IF.
\EQN
\BT\label{TH2}
Let $X(t),t\inr$ be a centered Gaussian process with continuous trajectories and stationary increments
satisfying {\bf AI-AII} {and $c>0$}. If
$$\varphi=\lim_{u\rw\IF}\frac{\sigma^2(u)}{u}\in [0,\IF],$$
{then 
for any $x\ge 0$}
\BQNY
\ehe{p_\IF(u,x)}\sim \B_{X_{\varphi}}(x)\sqrt{\frac{2A\pi}{B}} \frac{u}{M(u)\overleftarrow{\sigma}
 \left(\frac{\sqrt{2}\sigma^2(ut^*)}{u(1+ct^*)}\right)}\Psi(M(u)),
\EQNY
where
\BQN\label{theta}
X_{\varphi}(t)=\left\{
            \begin{array}{ll}
B_{\re{2}\alpha_0}(t), & \hbox{if } \varphi=0 \\
\frac{1+ct^*}{\sqrt{2}\varphi t^*}X\left(\overleftarrow{\sigma}\left(\frac{\sqrt{2}\varphi t^*}{1+ct^*}\right)t\right), & \hbox{if }\varphi\in (0,\IF)\\
B_{\re{2}\alpha_\IF}(t),&  \hbox{if } \varphi=\IF.
              \end{array}
            \right.
\EQN
\ET
\K{Application of Theorem \ref{TH2} to $X=B_1$ with comparison to (\ref{brownian}) leads to} the following corollary.
\BK For any $x\geq 0$,
\BQNY
\B_{B_1}(x)=(2+x)\Psi\left(\sqrt{\frac{x}{2}}\right)-\sqrt{\frac{x}{\pi}}e^{-\frac{x}{4}}.
\EQNY
\EK

\COM{\K{Now let us proceed to a complementary result that deals with the asymptotics of}
 $\tau_u(x)$, $x\geq 0$, as $u\to\infty$, defined by
\BQN\label{passengetime}\tau_u(x)=\inf\left\{t\geq 0: \frac{1}{\overleftarrow{\sigma}}\int_0^t \mathbb{I}_u(X(s)-cs)ds>x\right\}
\EQN
with the convention that $\inf\emptyset=\IF$. $\tau_u(x), x\geq 0$ are a family of  non-negative random variables satisfying that for fixed $u$, $\tau_u(x)$ is increasing with respect to $x$ for $x\in [0,\IF]$. Note that $\tau_u(0)$ is the first passage time of the process $X(t)-ct$  reaching level $u$.}

 \ehe{Next we analyze the asymptotic distribution of $\tau_u^*(x_1, x_2)$ defined  in
 \eqref{conditionaltime}, assuming that these random variables are defined on the same probability space.}
\BK\label{TH3} Under the assumptions of Theorem \ref{TH2}, for $0\leq x_1\leq x_2<\IF$,
the following convergence in distribution holds
\BQNY
\frac{\tau_u^*(x_1,x_2)-ut_u}{A(u)}\stackrel{d}{\rw}\mathcal{N}_{x_1,x_2}, \quad u\rw\IF,
\EQNY
where
$$ \pk{\mathcal{N}_{x_1,x_2}\leq y}=\frac{\B_{X_{\varphi}}(x_2)}{\B_{X_{\varphi}}(x_1)}\pk{\mathcal{N}\leq y}, \quad y\in \mathbb{R}, \quad \pk{\mathcal{N}_{x_1,x_2}=\IF}=1-\frac{\B_{X_{\varphi}}(x_2)}{\B_{X_{\varphi}}(x_1)},$$
with
$A(u)=\frac{\sigma(ut^*)}{c}\sqrt{\frac{\alpha_\IF}{1-\alpha_\IF}}$ and $\mathcal{N}$ an $ N(0,1)$  random variable.
\EK
\subsubsection{Finite time horizon.}
{In this subsection we consider the finite-time horizon case, i.e., we are interested in the asymptotics of $ p_T(u,x)$
 as $u\to \IF$, where $X$ has stationary increments$,  x\ge 0$, $T\in (0,\IF)$ {and $c\in\R$.}
\K{Due to the finiteness of $T$, we allow in this section $c\in \mathbb{R}$.}
Clearly,
\BQN\label{transfinite}
p_T(u,x)&=&\pk{v(u)\int_0^T\mathbb{I}_0\left(\frac{X(t)}{u+ct}-1\right)\td t>x}\nonumber\\
&=&\pk{v(u)\int_0^T\mathbb{I}_{\K{m(u)}}\left(\frac{X(t)}{u+ct}m(u)\right)\td t>x}.
\EQN
where
$$m(u)=\frac{u+cT}{\sigma(T)}.$$ 

We shall impose the following assumptions on $\sigma$.\\
\\
{\bf BI} $\sigma(0)=0$ and $\sigma\in C([0,T])$  with the first derivative $\dot{\sigma}(t)>0 , t\in (0,T]$.\\
{\bf BII} $\sigma$ is regularly varying at $0$ with index $\alpha_0\in (0,1]$.\\
\\
{We note that both fBm and introduced in Section \ref{s.sti} Gaussian integrated processes $\int_0^tZ(s)\td s$
with correlation of $Z$ such that $r(t)>0$ satisfy conditions {\bf BI-BII}
(note that $\sigma^2(t)=2\int_0^t\int_0^s r(u)\td u \td s,~t\geq 0$ and $\alpha_0=1$).}
%

Assumption {\bf BI} ensures that the first derivative of $\frac{\sigma(t)}{u+ct}$ is positive and further its maximizer
over $[0,T]$ is unique and \re{equals} $T$ for sufficiently large u. {Assumption} {\bf BII} gives the correlation structure of $\frac{X(t)}{u+ct}$
around time $t=T$ {(see  Lemma \ref{LemL2} ii)}.
Thus under the assumptions of {\bf BI}-{\bf BII}, $\frac{X(t)}{u+ct}m(u), t\in [0,T]$
is a centered continuous Gaussian process with the maximum of variance function
\kk{attained at $t=T$ and equal to $1$.} \\
 Set
$$ v(u)=\left\{\begin{array}{cc}
1/\overleftarrow{\sigma}\left(\frac{\sqrt{2}\sigma^2(T)}{u+cT}\right),& \lim_{t\rw 0}\frac{|t|}{\sigma^2(|t|)}\in [0,\IF)\\
(m(u))^2,& \lim_{t\rw 0}\frac{|t|}{\sigma^2(|t|)}=\IF.
\end{array}
\right.$$

\BT\label{TH4} Suppose that $X(t), t\in [0,T]$ is a centered continuous Gaussian process with stationary increments satisfying
{\bf BI-BII} \ehd{and let $x\ge 0${, $c\in\R$} be given.}\\
i) If $t=o(\sigma^2(t))$ as $t\rw 0$, then
\BQNY p_T(u,x)\sim \B_{B_{2\alpha_0}}(x)\frac{\sigma(T)}{\dot{\sigma}(T)}\frac{1}{(m(u))^2\overleftarrow{\sigma}\left(\frac{\sqrt{2}\sigma^2(T)}{u+cT}\right)}\Psi(m(u)).
\EQNY	
ii) If $ \lim_{t\rw 0} \sigma^2(t)/t= \theta \kk{\in(0,\infty)} $, then
\BQNY p_T(u,x)\sim \B_{B_{1}}^{\frac{2\sigma(T)\dot{\sigma}(T)}{\theta}|t|}(x)\Psi(m(u)).
\EQNY
iii) If
$\sigma^2(t)=o(t)$ as $t\rw 0$, then
\BQNY
p_T(u,x) \sim e^{-\frac{\dot{\sigma}(T)}{\sigma(T)}x}\Psi(m(u)).
\EQNY
\ET

With the convention that $\inf\emptyset=\IF$,
define  $\tau_{u,T}(x)$, $T>0,x\geq 0$ by
\BQN
\tau_{u,T}(x)=\inf\left\{t: v(u)\int_0^t \mathbb{I}_u(X(s)-cs)\td s>x, 0\leq t\leq T\right\}.
\EQN
Further, \K{let}
\BQN\label{Fintetime}\tau_{u,T}^*(x_1,x_2):=\tau_{u,T}(x_2)\Big {|}\tau_{u,T}(x_1)\leq T, \quad 0\leq x_1\leq x_2<\IF.
\EQN

\BK\label{TH5}
\kk{Suppose that $X(t), t\in [0,T]$ is a centered continuous Gaussian process with stationary increments satisfying {\bf BI-BII},
  \re{$0\leq x_1\leq x_2<\IF$, {$c\in\R$} are given}
and}
 $$\lim_{t\rw 0}\frac{|t|}{\sigma^2(|t|)}\in [0,\IF].$$
Then
$$\frac{\dot{\sigma}(T)}{\sigma^3(T)}u^2(T-\tau_{u,T}^*(x_1,x_2))\stackrel{d}{\rw}\mathcal{E}_{x_1,x_2}, \quad 0\leq x_1\leq x_2<\IF, \quad u\rw\IF,$$
where
$$\pk{\mathcal{E}_{x_1,x_2}>y}=\Gamma(x_1,x_2)e^{-y}, \quad \pk{\mathcal{E}_{x_1,x_2}=-\IF}=1-\Gamma(x_1,x_2)\geq 0,\quad \re{y\geq 0},$$
with
$$\Gamma(x_1,x_2)=\left\{\begin{array}{cc}
\frac{\B_{B_{2\alpha_0}}(x_2)}{\B_{B_{2\alpha_0}}(x_1)}& t=o(\sigma^2(t))\\
\frac{\B_{B_{1}}^{\frac{2\sigma(T)\dot{\sigma}(T)}{\theta}|t|}(x_2)}{\B_{B_{1}}^{\frac{2\sigma(T)\dot{\sigma}(T)}{\theta}|t|}(x_1)} & \sigma^2(t)\sim \theta t\\
e^{\frac{\dot{\sigma}(T)}{\sigma(T)}(x_1-x_2)}& \sigma^2(t)=o(t)
\end{array}\right., \quad 0\leq x_1\leq x_2<\IF.$$
\EK
\subsection{Self-similar Gaussian processes} {In this subsection we apply our findings
to the class of self-similar Gaussian processes with drift.
We focus on the exact asymptotic probabilities of sojourn time,
without specifying analogs of Corollaries \ref{TH3} and \ref{TH5},
since they essentially lead to the same type of results as given in Section \ref{s.sta}.
}

\K{Suppose that $X$ is a centered self-similar Gaussian process
\kk{with self-similarity index $H\in (0,1)$}, i.e.}
\BQN\label{selfsimilar}
\{X(bt), t\geq 0\}\stackrel{d}{=} \{b^H X(t), t\geq 0\},
\EQN
where $\stackrel{d}{=}$ means the equality of finite dimensional distributions.

\K{Equality (\ref{selfsimilar}) implies} that $\sigma^2(t)=Var(X(t))=Var(X(1))t^{2H}$.
Without loss of generality, \kk{in what follows} we assume that $Var(X(1))=1$ and \kk{hence} $\sigma^2(t)=t^{2H}$.\\
Moreover, we assume that

{\bf S} There exist a function $\rho$ which is regularly varying at $0$ with index $\alpha\in (0,2]$, $\rho(0)=0$ and
$\rho(t)>0, t>0$, and $t_0\in [0,T]$ such that
\BQN\label{selfcor}
\lim_{\epsilon\rw 0}\sup_{s\neq t,  \blu{|s-t_0|<\epsilon}, |t-t_0|<\epsilon, s,t\in [0,T]}\left|\frac{1-Corr(X(s), X(t))}{\rho(|t-s|)}-1\right|=0.
\EQN

{Condition {\bf S} is satisfied by such classes of self-similar Gaussian processes as
 fBms, bi-fractional Brownian motions, sub-fractional Brownian motions or generalized fractional Brownian
 motions; see, e.g.,
 \cite{Mishura18}, \cite{BG2004}, \cite{Zili18} or \cite{hj2015}.
}
\subsubsection{Infinite-time horizon}
{Suppose that $T=\IF$ {and $c>0$.} Then,}
by self-similarity of $X$, we have
$$\pk{v(u)\int_0^\IF \mathbb{I}_u(X(t)-ct)dt>x}=
\pk{uv(u)\int_0^\IF \mathbb{I}_{\K{u^{1-H}}}\left(\frac{X(t)}{1+ct}\right)dt>x}.$$
Note that the maximizer  of
$\sqrt{Var\left(\frac{X(t)}{1+ct}\right)}=\frac{t^H}{1+ct}$ is unique and \K{equals}
$\frac{H}{c(1-H)}$. Further, referring to \cite{HP99},
\BQN\label{selfvariance}
\widehat{A}\frac{t^H}{1+ct}=1-\frac{\widehat{B}}{2\widehat{A}}\left(t-\frac{H}{c(1-H)}\right)^2(1+o(1)), \quad t\rw \frac{H}{c(1-H)},
\EQN
with
$$\widehat{A}=\left(\frac{{H}}{c(1-{H})}\right)^{-H}\frac{1}{1-H}, \quad
\widehat{B}=\left(\frac{H}{c(1-H)}\right)^{-H-2}H.$$
We arrive at the following result.
\BT\label{selfsimilar0} Let $X(t)$ be a centered self-similar Gaussian process
\kk{with self-similarity index $H\in(0,1)$} satisfying
(\ref{selfsimilar}) and {\bf S} with $t_0=\frac{H}{c(1-H)}$. \K{Suppose} that
$$\lim_{t\rw 0}\frac{t^2}{\rho(t)}=\gamma\in [0,\IF]$$
\ehd{and let $x\ge 0${, $c>0$} be given}. \\
\K{i)} If $\gamma=0$, then
\BQNY
\pk{\frac{1}{u\overleftarrow{\rho}((\widehat{A}u^{1-H})^{-2})}\int_0^\IF \mathbb{I}_u(X(t)-ct)dt>x}\sim \B_{B_\alpha}(x)\sqrt{\frac{2\widehat{A}\pi}{\widehat{B}}} \frac{1}{\overleftarrow{\rho}((\widehat{A}u^{1-H})^{-2})\widehat{A}u^{1-H}}\Psi(\widehat{A}u^{1-H}).
\EQNY
\K{ii)} If $\gamma\in (0,\IF)$, then
\BQNY
\pk{\frac{1}{u\overleftarrow{\rho}((\widehat{A}u^{1-H})^{-2})}\int_0^\IF \mathbb{I}_u(X(t)-ct)dt>x}\sim \sqrt{\frac{2\widehat{A}+\widehat{B}\gamma}{\widehat{B}\gamma}}e^{-\frac{2\widehat{A}+\widehat{B}\gamma}{8\widehat{A}}x^2}\Psi(\widehat{A}u^{1-H}).
\EQNY
\K{iii)} If $\gamma=\IF$, then
\BQNY
\pk{\blu{u^{-H}}\int_0^\IF \mathbb{I}_u(X(t)-ct)dt>x}\sim e^{-\frac{\widehat{A}\widehat{B}x^2}{8}}\Psi(\widehat{A}u^{1-H}).
\EQNY
\ET
\COM{Let
$$v(u)=\left\{\begin{array}{cc}
\frac{1}{u\overleftarrow{\rho}((\widehat{A}u^{1-H})^{-2})}& \gamma\in [0,\IF)\\
u^{-H}& \gamma=\IF.
\end{array}
\right.$$

\K{Finally, we consider $\tau_u^*(x_1,x_2)$ defined in (\ref{conditionaltime}), where
we substitute} $ \frac{1}{\overleftarrow{\sigma}}$ \blu{for} $v(u)$ in (\ref{passengetime}), \K{that is we} have that for $x\geq 0$
$$\tau_u(x):=\inf\left\{t\geq 0: v(u)\int_0^t \mathbb{I}_u(X(s)-cs)ds>x\right\}.$$
\BK\label{selfsimilarad}
\K{Suppose that}
conditions of Theorem \ref{selfsimilar0}, with $\tau_u^*(x_1,x_1)$ defined in (\ref{conditionaltime}), \K{are satisfied}
for $0\leq x_1\leq x_2<\IF$. \\
\K{i)} If $\gamma=0$, then
\BQNY
\frac{\sqrt{\widehat{A}\widehat{B}}(\tau_u^*(x_1,x_2)-ut_0)}{u^{H}}\stackrel{d}{\rw}\mathcal{N}_{x_1,x_2}, \quad u\rw\IF,
\EQNY
where
$$ \pk{\mathcal{N}_{x_1,x_2}\leq y}=\frac{\B_{B_\alpha}(x_2)}{\B_{B_\alpha}(x_1)}\pk{\mathcal{N}\leq y}, \quad y\in \mathbb{R}, \quad \pk{\mathcal{N}_{x_1,x_2}=\IF}=1-\frac{\B_{B_\alpha}(x_2)}{\B_{B_\alpha}(x_1)},$$
with
 $\mathcal{N}$ an $ N(0,1)$
 random variable. \\
\K{ii)} If $\gamma\in (0,\IF)$, then for any $y\in (-\IF, \IF]$
 \BQNY
 \limit{u}\pk{\frac{(\tau_u^*(x_1,x_2)-ut_0)}{(\widehat{A})^{-1}\sqrt{\gamma}u^{H}}\leq y} \ehe{=} e^{\frac{2\widehat{A}+\widehat{B}\gamma}{8\widehat{A}}x_1^2}F_{\frac{\widehat{B}\gamma}{2\widehat{A}}}(x_2,y),
 \EQNY
 where $F$ is defined in (\ref{Pit1}).\\
\K{iii)} If $\gamma=\IF$, then for any $y\in \mathbb{R}$
\BQNY
\limit{u}\pk{\sqrt{\frac{\widehat{A}\widehat{B}}{2}}\frac{(\tau_u^*(x_1,x_2)-ut_0)}{u^{H}}\leq y}=\left\{\begin{array}{cc}
e^{\frac{\widehat{A}\widehat{B}x_1^2}{8}-(x_2\sqrt{\frac{\widehat{A}\widehat{B}}{2}}-y)^2},& y\leq x_2\sqrt{\frac{\widehat{A}\widehat{B}}{8}}\\
e^{\frac{\widehat{A}\widehat{B}(x_1^2-x_2^2)}{8}}, &y> x_2\sqrt{\frac{\widehat{A}\widehat{B}}{8}}
\end{array}
\right..
\EQNY
\EK}
\subsubsection{Finite-time horizon} Let $T\in (0,\IF)$. We arrive at the following result.

\BT\label{selfsimilar1}
\kk{Let $X(t)$ be a centered self-similar Gaussian process
\kk{with self-similarity index $H\in(0,1)$} satisfying
(\ref{selfsimilar}) and \K{{\bf S}} with $t_0=T$.}
\ehd{For given $x\ge 0$ { and $c\in\R$,} let}  $ p_T(u,x)$ be defined in (\ref{PT}) with $v(u)=\overleftarrow{\rho}\left(\frac{T^{2H}}{(u+cT)^2}\right)$ if $\lim_{t\rw 0}\frac{|t|}{\rho(|t|)}\in [0,\IF)$ and $v(u)=\frac{(u+cT)^2}{T^{2H}}$ if $\lim_{t\rw 0}\frac{|t|}{\rho(|t|)}=\IF$. \\
i) If $t=o(\rho(t))$ as $t\rw 0$, then
\bqny{ p_T(u,x)\sim \B_{B_{\alpha}}(x)\frac{T^{2H+1-2H/\alpha}}{H}\frac{1}{u^2\overleftarrow{\rho}(u^{-2})}\Psi\left(\frac{u+cT}{T^H}\right).
}	
ii) If $ \lim_{t\rw 0} \rho(t)/t= \theta $, then
\bqny{ p_T(u,x)\sim \B_{B_{1}}^{\frac{H}{T\theta}|t|}(x)\Psi\left(\frac{u+cT}{T^H}\right).
}
iii) If
$\rho(t)=o(t)$ as $t\rw 0$, then
\bqny{p_T(u,x) \sim e^{-\frac{H}{T}x}\Psi\left(\frac{u+cT}{T^H}\right).
}
\ET
\COM{\BK\label{TH50} \ehe{If $X$ is as in \netheo{selfsimilar1}} and
 $$\lim_{t\rw 0}\frac{|t|}{\rho(|t|)}\in [0,\IF],$$
then  the following convergence in distribution \K{holds}
$$\frac{\dot{\sigma}(T)}{\sigma^3(T)}u^2(T-\tau_{u,T}^*(x_1,x_2))\stackrel{d}{\rw}\mathcal{E}_{x_1,x_2}, \quad 0\leq x_1\leq x_2<\IF, \quad u\rw\IF,$$
where
$$\pk{\mathcal{E}_{x_1,x_2}>y}=\Gamma(x_1,x_2)e^{-y}, \quad \pk{\mathcal{E}_{x_1,x_2}=\IF}=1-\Gamma(x_1,x_2)\geq 0,$$
with
$$\Gamma(x_1,x_2)=\left\{\begin{array}{cc}
\frac{\B_{B_{\alpha}}(x_2)}{\B_{B_{\alpha}}(x_1)}& t=o(\rho(t))\\
\frac{\B_{B_{1}}^{\frac{H}{T\theta}|t|}(x_2)}{\B_{B_{1}}^{\frac{H}{T\theta}|t|}(x_1)} & \rho(t)\sim \theta t\\
e^{\frac{H}{T}(x_1-x_2)}& \rho(t)=o(t)
\end{array}\right., \quad 0\leq x_1\leq x_2<\IF.$$
\EK}

\section{Technical Lemmas}\label{s.lemmas}
We first present \ehe{a modification} of Theorem 5.1 \K{in} \cite{KEZX17}, which is crucial for the proofs  below. Then we present two lemmas related with the local behavior of the variance and correlation function of the investigated Gaussian process.

Let $\{\K{\xi_{u,k}}(t), t\in E,k\in K_u\}$, with  $K_u$ an index set and $E=[a,b]$  with $a\leq 0\leq b$ be a family of centered continuous Gaussian random \re{processes} with  variance function $\sigma^2_{\xi_{u,k}}$.
We impose the following  assumptions:\\
\\
{\bf C0}:
 Let $\{\K{g}_{k}(u),k\in K_u\}$ be a sequence of deterministic functions of $u$ satisfying
     \BQNY
     \lim_{u\to\IF}\inf_{k\in K_u}g_{k}(u)=\IF.
     \EQNY
{\bf C1}: $\sigma_{\xi_{u,k}}(0)=1$ for all large $u$ and there exists a \re{continuous} function $h$ \re{on $E$} such that
\BQNY\label{assump-cova-field}
\lim_{u\to\IF}\sup_{t\in E,k\in K_u}\ABs{g_k^2(u)\LT( 1-\sigma_{\xi_{u,k}}(t) \RT) - h(t) } =0.
\EQNY	
{\bf C2}: \K{There exists a centered Gaussian process with continuous trajectories and stationary increments $\zeta(t),t\inr$,} satisfying {\bf AI-AII}, \re{and}  \BQNY
\lim_{u\rw\IF}\sup_{k\in K_u}\sup_{s\neq t, s,t\in E}\left|\frac{g_k^2(u)\left(1-Corr\left(\xi_{u,k}(s), \xi_{u,k}(t)\right)\right)}{\frac{\sigma_\zeta^2(\upsilon(u)\Delta(u)|t-s|)}{\sigma_\zeta^2(\Delta(u))}}-1\right|=0,
\EQNY
where
 $\Delta(u),\upsilon(u),u>0$ are positive functions such that
\begin{equation}\label{DDelta}\lim_{u\rw\IF}\Delta(u)=\varphi\in [0,\IF], \quad \lim_{u\rw\IF}\upsilon(u)=\upsilon\in [0,\IF).
\end{equation}
\COM{{\bf C2}: These exists a centered continuous Gaussian process $\zeta(t), t\in\mathbb{R}$ with $\zeta(0)=0$ a.e., such that for any $s,t\in E$
\BQNY\label{assump-vari-field}
  \lim_{u\to\IF}\sup_{k\in K_u}\ABs{g_k^2(u)\big(Var(\overline{\xi}_{u,k}(t)-\overline{\xi}_{u,k}(s))\big) - 2Var(\zeta(t)-\zeta(s)) } =0.
\EQNY
{\bf C3}: These exist positive constants $C, \nu, u_0$ such that
\BQNY\label{assump-holder-field}
  \sup_{k\in K_u}  g_k^2(u)\E[\Big]{ (\overline{\xi}_{u,k}(t)-\overline{\xi}_{u,k}(s))^2} \leq C|t-s|^\nu
  \EQNY}
\K{Let}, \ehd{for given $x\ge 0$} and continuous function $h$
\BQNY
\B^{h}_{\zeta_{\varphi}}(x,E)= \int_\R  \pk{ \int_{E} \mathbb{I}_0\big(\sqrt{2}\zeta_{\varphi}(t)-\sigma_{\zeta_{\varphi}}^2(t)-h(t) + z \big)\td t> x } e^{-z} \td z,
\EQNY
with
$$\zeta_{\varphi}(t)=\left\{\begin{array}{cc}
B_{2\alpha_0}(\upsilon t)& \varphi=0\\
\frac{\zeta(\upsilon\varphi t)}{\sigma_\zeta(\varphi)}& \varphi\in (0,\IF)\\
B_{2\alpha_\IF}(\upsilon t)& \varphi=\IF
\end{array}
\right..$$
In the special case   $\upsilon=0$,  $\zeta_{\varphi}(t)\equiv 0$. Note that by Borell-TIS Inequality \cite{GennaBorell, AdlerTaylor}, we have
\BQNY
\B^{h}_{\zeta_{\varphi}}(x,E)&\leq& \B^{h}_{\zeta_{\varphi}}(0,E)\\
&=&\int_\R\pk{ \sup_{t\in E} \left(\sqrt{2}\zeta_{\varphi}(t)-\sigma_{\zeta_{\varphi}}^2(t)-h(t)\right)>z } e^{z} \td z\\
&\leq & e^a+\int_a^\IF e^{-\frac{(z-a)^2}{2b}}e^z\td z<\IF\,,
\EQNY
with
$$a=\mathbb{E}\left(\sup_{t\in E}\left(\sqrt{2}\zeta_{\varphi}(t)-\sigma_{\zeta_{\varphi}}^2(t)-h(t)\right)\right)<\IF,\quad  b=\sup_{t\in E} Var(\sqrt{2}\zeta_{\varphi}(t))<\IF.$$
\BEL\label{Pickands}
Let $\{\xi_{u,k}(t), t\in E,k\in K_u\}$ with $E=[a,b]$ be a family of centered continuous Gaussian \K{processes} satisfying {\bf C1-C2}. If
\K{ $g_k(u),k\in K_u$ satisfies} {\bf C0} and
 for any $x\geq0$
$$ \pk{ \int_{E} \mathbb{I}_0\LT(\xi_{u,k}(t)-g_k(u)\RT)\td t >x} >0,\quad \forall\,k\in K_u,$$
\K{then}
\BQN\label{con-uni-con}
\lim_{u\to\IF}\sup_{k\in K_u} \left| \frac{\pk{ \int_{E} \mathbb{I}_0\LT(\xi_{u,k}(t)-g_k(u)\RT) \td t>x }  } {\Psi(g_k(u))} -\B^{h}_{\zeta_{\varphi}}(x,E) \right| = 0
\EQN
holds for all $x\geq 0$.
Additionally, $\B^{h}_{\zeta_{\varphi}}(x,E)$ is a continuous function over $ [0, mes(E))$.

\EEL
\BRM i) In the special {case $\upsilon=0$, for all $x\geq 0$, we have}
\BQN\label{con-uni-con1}
\lim_{u\to\IF}\sup_{k\in K_u} \left| \frac{\pk{ \int_{E} \mathbb{I}_0\LT(\xi_{u,k}(t)-g_k(u)\RT) \td t>x }  } {\Psi(g_k(u))} -\B^{h}_{0}(x,E) \right| = 0.
\EQN
ii) There exists a non-negative  function $\rho$ which is a regularly varying function at $0$ with index $2\alpha_0\in (0,2]$  such that
 \BQN\label{replace}
\lim_{u\rw\IF}\sup_{k\in K_u}\sup_{s\neq t, s,t\in E}\left|\frac{g_k^2(u)\left(1-Corr\left(\xi_{u,k}(s), \xi_{u,k}(t)\right)\right)}{\frac{\rho(\upsilon(u)\Delta(u)|t-s|)}{\rho(\Delta(u))}}-1\right|=0,
\EQN
with $\Delta(u),\upsilon(u),u>0$ positive functions {satisfying (\ref{DDelta}) for $\varphi=0$.}
 If {\bf C2} is replaced by (\ref{replace}), then Lemma \ref{Pickands} still holds with $\zeta_{\varphi}(t)$ replaced by $B_{2\alpha_0}(\upsilon t)$.
\ERM
\prooflem{Pickands}. In order to prove this lemma it suffices to check the conditions of \blu{Theorem 5.1} in \cite{KEZX17}. That is we have to prove that (with $\overline{\xi}_{u,k}(t)=\frac{\xi_{u,k}(t)}{\sigma_{\xi_{u,k}}(t)}$)

i)
     $
     \lim_{u\to\IF}\inf_{k\in K_u}g_k(u)=\IF.
     $\\
ii) $\sigma_{\xi_{u,k}}(0)=1$ for all large $u$ and any $k\in K_u$, and
there exists some bounded continuous function $h$ on $E$ such that
\BQNY\label{assump-cova-field}
\lim_{u\to\IF}\sup_{t\in E,k\in K_u}\abs{g_k^2(u)\LT( 1- \E{ \xi_{u,k}(t)\xi_{u,k}(0) } \RT) -\sigma_{\zeta_{\varphi}}^2(t)- h(t) } =0.
\EQNY	

iii) For any $s,t\in E$,
\BQNY\label{assump-vari-field}
  \lim_{u\to\IF}\sup_{ k\in K_u}\abs{g_k^2(u)\big(Var(\overline{\xi}_{u,k}(t)-\overline{\xi}_{u,k}(s))\big) - 2Var(\zeta_{\varphi}( t)-\zeta_{\varphi}( s)) } =0.
\EQNY
iv) These exist positive constants $C, \nu, u_0$ such that
\BQNY\label{assump-holder-field}
  \sup_{k\in K_u}  g_k^2(u)\E{\overline{\xi}_{u,k}(t)-\overline{\xi}_{u,k}(s)}^2 \leq C \norm{s-t}^\nu
  \EQNY
  holds for all $s,t\in E, u\geq u_0$.\\
It follows that i)  holds straightforwardly from {\bf C0}. Next we verify ii), iii) and iv).\\
 We first prove iii).
Uniform convergence theorem (see, e.g., \cite{BI1989}) gives that, for all $\upsilon\in [0,\IF)$ and $\varphi\in \{0,\IF\}$,
\BQN\label{unicon}\lim_{u\rw\IF}\sup_{s,t\in E}\left|\frac{\sigma_{\zeta}^2(\upsilon(u)\Delta(u)|t-s|)}{\sigma_{\zeta}^2(\Delta(u))}-Var\left(\zeta_{\varphi}( t)-\zeta_{\varphi}(\re{s})\right)\right|=0.
\EQN
For $\upsilon\in [0,\IF)$ and $\varphi\in (0,\IF)$, the above limit is still valid by the continuity of $\sigma_{\zeta_\varphi}^2$.
Hence it follows from  {\bf C2} that
\BQN\label{a1}
\lim_{u\rw\IF}\sup_{k\in K_u}\sup_{ s,t\in E}\left|g_k^2(u)Var\left(\overline{\xi}_{u,k}(s)-\overline{\xi}_{u,k}(t)\right)-2Var\left(\zeta_{\varphi}(s)-\zeta_{\varphi}(t)\right)\right|=0.
\EQN
This confirms that iii) is satisfied.\\
We next focus on iv). Let $f(t)=\frac{\sigma_\zeta^2(t)}{t^\lambda}$ with $0<\lambda<2\min (\alpha_0, \alpha_\IF)$. Then $f$ is a regularly varying function at $0$ and $\IF$ with index $2\alpha_0-\lambda>0$ and $2\alpha_\IF-\lambda>0$, respectively.
By {\bf C2}, we have that, for $u$ sufficiently large,
\BQNY
\sup_{k\in K_u}  g_k^2(u)\E{ (\overline{\xi}_{u,k}(t)-\overline{\xi}_{u,k}(s))^2}&\leq& 4\frac{\sigma_\zeta^2(\upsilon(u)\Delta(u)|t-s|)}{\sigma_\zeta^2(\Delta(u))}\\
&\leq& 4(\upsilon+1)^\lambda\frac{f(\Delta(u)\upsilon(u)|t-s|)}{f(\Delta(u))}|t-s|^\lambda , \quad s,t\in E.
\EQNY
 Uniform convergence theorem yields that for $0<\lambda<2\min (\alpha_0, \alpha_\IF)$ and $\varphi\in \{0,\IF\}$, $\upsilon\in [0,\IF)$
\BQNY
\re{\lim_{u\to\IF}}\sup_{s\neq t, s,t\in E}\left|\frac{f(\Delta(u)\upsilon(u)|t-s|)}{f(\Delta(u))}-(\upsilon|t-s|)^{(2\alpha_0-\lambda)
I_{\{\varphi=0\}}+(2\alpha_\IF-\lambda)I_{\{\varphi=\IF\}}}\right|=0,
\EQNY
which {implies} that for $u$ large enough
$$\sup_{s\neq t, s,t\in E}\frac{f(\Delta(u)\upsilon(u)|t-s|)}{f(\Delta(u))}<\mathbb{Q}.$$
with $\mathbb{Q}$ a positive constant.
For $\varphi\in (0,\IF)$ and $\upsilon\in (0,\IF)$ the above inequality follows straightforwardly by noting the continuity of $f$.
Hence for $u$ sufficiently large, for all $\varphi \in [0,\IF]$ and $\upsilon\in [0,\IF)$
\BQN\label{a2}
\sup_{k\in K_u}  g_k^2(u)\E{ (\overline{\xi}_{u,k}(t)-\overline{\xi}_{u,k}(s))^2}\leq \mathbb{Q}|t-s|^\lambda , \quad s,t\in E,
\EQN
with $0<\lambda<2\min (\alpha_0, \alpha_\IF)$. This implies that iv) hold.\\
Finally, we prove ii).   Notice that
$$ 1- \E{ \xi_{u,k}(t)\xi_{u,k}(0) }=(1-\sigma_{\xi_{u,k}}(t))+(1-Corr(\xi_{u,k}(t),\xi_{u,k}(0)))-(1-\sigma_{\xi_{u,k}}(t))(1-Corr(\xi_{u,k}(t),\xi_{u,k}(0))).$$
In view of (\ref{a1}),
$$\lim_{u\to\IF}\sup_{t\in E,k\in K_u}\abs{g_k^2(u)(1-Corr(\xi_{u,k}(t),\xi_{u,k}(0)))-\sigma_{\zeta_{\varphi}}^2(t)}=0,$$
which combined with {\bf C1} leads to
$$\lim_{u\to\IF}\sup_{t\in E,k\in K_u}\abs{g_k^2(u)(1-Corr(\xi_{u,k}(t),\xi_{u,k}(0)))(1-\sigma_{\xi_{u,k}}(t))}=0.$$
Hence in view of  {\bf C1} and the above limits,
$$\lim_{u\to\IF}\sup_{t\in E,k\in K_u}\abs{g_k^2(u)(1- \E{ \xi_{u,k}(t)\xi_{u,k}(0) })-\sigma_{\zeta_{\varphi}}^2(t)-h(t)}=0.$$
This confirms that ii) holds.
Thus, {applying \blu{Theorem 5.1} in \cite{KEZX17},} the claim is established for all the continuity points of $\mathcal{B}^{h}_{\zeta}( x, E)$.\\
{\it \underline{Continuity of $\mathcal{B}^{h}_{\zeta}( x, E)$}}.
Next we show that $\mathcal{B}^{h}_{\zeta}( x, E)$ is continuous over $ [0, mes(E))$.
{Since $\mathcal{B}^{h}_{\zeta}( x, E)$ is right-continuous at $0$, then} we are left with the continuity of  $\mathcal{B}^{h}_{\zeta}( x, E)$  over $(0, mes(E))$.    The claimed continuity at $x\in (0, mes(E))$ follows if
$$\int_\R  \pk{ \int_{ E } \mathbb{I}\big(\sqrt{2}\zeta_{\varphi}(t)-\sigma_{\zeta_{\varphi}}^2(t)-h(t) + s >0\big) \td t =x } e^{-s} \td s=0.$$
Assume that the probability space $(C(E), \mathcal{F}, \mathbb{P}^*)$ with $E=[a,b]$  is induced by the process  $\{\sqrt{2}\zeta_{\varphi}(t)-\sigma_{\zeta_{\varphi}}^2(t)-h(t) , t\in E\}$ with $C(E)$ denoting the collection of all continuous functions over $E$ and equipped with sup-normal, and $\mathcal{F}$ being the Borel $\sigma$-field \re{on} $C(E)$. Thus in order to complete the proof, it suffices to prove that for any $x\in (0, mes(E))$
$$\int_\R  \mathbb{P}^*\left({ \int_{ E } \mathbb{I}\big(w(t) + s >0\big) \td t =x }\right) e^{-s} \td s=0,$$
where $w\in C(E)$.  For any $x\in (0, mes(E))$, let
$$A_s=\left\{ w\in C(E)\re{:} \int_{ E} \mathbb{I}\big(w(t)+ s >0\big) \td t =x \right\}, \quad s\in\mathbb{R}.$$
If $\int_{ E}\mathbb{I}\big(w(t)+ s >0\big) \td t =x$ with $x\in (0, mes(E))$, then $\inf_{t\in E} w(s)\leq-s<\sup_{t\in E} w(t)$. By the continuity of $w$, for $s<s'$,
$$\int_{ E}\mathbb{I}\big(w(t)+ s' >0\big) \td t\in (x, mes(E)].$$
This implies that
\BQNY
A_s\cap A_{s'}=\emptyset, \quad s\neq s', s,s'\in \mathbb{R}.
\EQNY
 Since  $A_s, s\in\mathbb{R}$ are  measurable sets and
{$$\sup_{\Lambda\subset \mathbb{R}, \#\Lambda<\infty}\sum_{s\in\Lambda}\mathbb{P}^*\left(A_s\right)\leq 1,$$}
where $\#{\Lambda}$ means the cardinality of set ${\Lambda}$, then
$$\{s: s\in \mathbb{R} \quad \text{such that} \quad \mathbb{P}^*\left({A_s}\right)>0\}$$
is a countable set, which indicates that
$$\int_\R  \mathbb{P}^*\left({A_s}\right) e^{-s} \td s=0.$$
Hence $\mathcal{B}^{h}_{\zeta_\varphi}( x, E)$ is continuous over $(0, mes(E))$.  This completes the proof.\QED

\bigskip

Let
$$\sigma_u^*(t)=Var^{1/2}\left(\frac{X(ut)}{u(1+ct)}M(u)\right), t\geq 0 \quad r_u(s,t)=Corr\left(\frac{X(us)}{u(1+cs)}, \frac{X(ut)}{u(1+ct)}\right), \quad s,t> 0.$$
Assume  that $\delta_u>0$ with $\lim_{u\rw\IF}\delta_u=0$.
\K{The following lemma is due to \cite{KrzysPeng2015}.}
\BEL \label{LemL1}
i) Suppose that {\bf AI} is satisfied and let $t_u=\arg\inf_{t>0}\frac{u(1+ct)}{\sigma(ut)}$.
 If $u$ is large enough then $t_u$ is unique and  $\limit{u} t_u= t^*=\frac{\alpha_\IF}{c(1-\alpha_\IF)}$. Moreover,
\BQN \label{L1}
\lim_{u\rw\IF}\sup_{t\in (t_u-\delta_u, t_u+\delta_u)\setminus \{t_u\}}\left|\frac{1-\sigma_u^*(t)}{\frac{B}{2A}(t-t_u)^2}-1\right|=0,
\EQN
 with $A, B$ defined in (\ref{AB}).\\
ii) If {\bf AI, AII} hold,  then
 \BQN\label{L2}
\lim_{u\rw \IF}\sup_{s\neq t, s,t \in (t_u-\delta_u, t_u+\delta_u)}\left|\frac{1-r_u(s,t)}{\frac{\sigma^2(u|s-t|)}{\re{2\sigma^2(ut^*)}}}-1\right|=0.
\EQN
\EEL
Let next
$$\widehat{\sigma}_u(t)=Var^{1/2}\left(\frac{X(t)}{u+ct}m(u)\right), t\in [0,T], \quad \widehat{r}_u(s,t)=Corr\left(\frac{X(s)}{u+cs}, \frac{X(t)}{u+ct}\right), s,t\in [0,T].$$
\BEL \label{LemL2}  i) If {\bf BI} holds, then
\bqn{\label{L3}
	\lim_{u\rw\IF}\sup_{t\in (T-\delta_u, T)}\left|\frac{1-\widehat{\sigma}_u(t)}{|T-t|}-\frac{\dot{\sigma}(T)}{\sigma(T)}\right|=0.
}	
ii) If {\bf BII} holds and $t=o(\sigma(t))$ as $t\rw 0$,  then
 \BQN
 \label{L4}
\lim_{u\rw \IF}\sup_{s\neq t, s,t \in (T-\delta_u, T)}\left|\frac{1-\widehat{r}_u(s,t)}{\frac{\sigma^2(|s-t|)}{\re{2\sigma^2(T)}}}-1\right|=0.
\EQN
\EEL

\section{Proofs}\label{s.proofs}
Hereafter, denote by $\mathbb{Q}$, $\mathbb{Q}_i, i=1,2,3,\dots$ positive constants that may differ from line to line. The equivalence $f(u,S)\sim h(u)$ as  $u\rw\IF, S\rw\IF$ means that $\lim_{S\rw\IF}\lim_{u\rw\IF}\frac{f(u,S)}{h(u)}=1$. Moreover, for any non-constant random variable $X$, denote by $\overline{X}\coloneqq\frac{X}{\sqrt{Var(X)}}$.\\

 \prooftheo{ThZ}
By the lack of upward jumps and using the strong Markov property, for any $x, u\ge0$ we have
$$
\pk{\int_0^\IF \mathbb{I}_u(X(t)-ct)\td t >x} =
\pk{\int_0^\IF \mathbb{I}_0(X(t)-ct)\td t >x}\pk{\sup_{t\ge0}(X(t)-ct)\ge u}
$$
and using e.g. \cite{MandjesKrzys}[Theorem 3.3] we get $\pk{\sup_{t\ge0}(X(t)-ct)\ge u}=e^{-\alpha u}$, where
$\alpha$ is the positive solution to the equation
$\psi(\alpha)=c\alpha$
which by the strict convexity of $\psi$ and the assumption $\psi'(0+)< c$ exists and is unique.
This completes the proof.
\QED

\prooftheo{TH1} We consider each case i)-iii) separately. The following notation
\K{is valid} for \kk{all} three cases.\\
 Put for $u>0, x\ge 0$ \K{and $S>0$}
$$
\pi(u)=\pk{\int_{E(u)} \mathbb{I}_0(Z_u(t)-n(u))\td t>x},$$
$$\overline{Z}_{u,k}(t)=\overline{Z}_u(kS+t), \quad I_k=[kS, (k+1)S], \quad N_1(u)=\left[\frac{a_1(u)}{S}\right]-1, \quad N_2(u)=\left[\frac{a_2(u)}{S}\right]+1,$$
and for some $\epsilon \in (0,1)$ set
$$  n_{u,k}^-=n(u)\left(1+(1-\epsilon)\inf_{t\in I_k}w(g(u)|t|)\right),$$ \quad $$n_{u,k}^+=n(u)\left(1+(1+\epsilon)\sup_{t\in I_k}w(g(u)|t|)\right).$$
Note that
\BQNY
\pi(u)&\leq &\pk{\sum_{k=N_1(u)-2}^{N_2(u)+2}\int_{I_k}\mathbb{I}_0(Z_u(t)-n(u))\td t>x}\\
&\leq&\pk{\exists  N_1(u)-2\leq  k\leq N_2(u)+2 \quad \text{such that} \quad \int_{I_k}\mathbb{I}_0(Z_u(t)-n(u))\td t>x}\\
&&+ \pk{\exists  N_1(u)-2\leq  k,l\leq N_2(u)+2 , k\neq l\quad \text{such that} \quad \int_{I_k}\mathbb{I}_0(Z_u(t)-n(u))\td t>0, \int_{I_l}\mathbb{I}_0(Z_u(t)-n(u))\td t>0},
\EQNY
and
\BQNY
\pi(u)&\geq &\pk{\sum_{k=N_1(u)+2}^{N_2(u)-2}\int_{I_k}\mathbb{I}_0(Z_u(t)-n(u))\td t>x}\\
&\geq&\pk{\exists  N_1(u)+2\leq  k\leq N_2(u)-2 \quad \text{such that} \quad \int_{I_k}\mathbb{I}_0(Z_u(t)-n(u))\td t>x}.
\EQNY
Hence in view of (\ref{var}) and using Bonferroni inequality, we have
\BQN\label{Bonfer}
\pi^+(u,S)-\Sigma_1(u)-\Sigma_2(u)\leq \pi(u)\leq \pi^{-}(u,S)+\Sigma_1(u)+\Sigma_2(u),
\EQN
where
\BQNY\label{sigma1}
\pi^{\pm}(u,S)&=&\sum_{k=N_1(u)\pm 2}^{N_2(u)\mp 2}\pk{\int_{I_{\re{0}}}\mathbb{I}_0\left(\overline{Z}_{u,k}(t)-n_{u,k}^{\pm}\right)\td t>x},\nonumber\\
\Sigma_1(u)&=&\sum_{N_1(u)-2\leq k\leq N_2(u)+1}\pk{\sup_{t\in I_k}\overline{Z}_{u}(t)>n_{u,k}^-, \sup_{t\in I_{k+1}}\overline{Z}_{u}(t)>n_{u,\re{k+1}}^-},\\
\Sigma_2(u)&=&\sum_{N_1(u)-2\leq k<k+1<l\leq N_2(u)+1}\pk{\sup_{t\in I_k}\overline{Z}_{u}(t)>n^-_{u,k}, \sup_{t\in I_{l}}\overline{Z}_{u}(t)>n^-_{u,l}}.\nonumber
\EQNY

 { \underline{$\diamond$ \it Case i)}}
 The idea of the proof is to  divide $E(u)$ {into} a large number of tiny intervals \K{for each of which we}
 give a uniform  exact asymptotics of sojourn times.
For notational  simplicity define
$$\Theta(u)\coloneqq\beta^{-1}\int_{y_1}^{y_2} |t|^{1/\beta-1}e^{-|t|}\td t\frac{\inv{w}(n^{-2}(u))}{g(u)}\Psi(n(u)).$$ Without loss of generality, we assume that $x_1>0, \lim_{u\rw\IF} a_1(u)=\IF$ and $x_2>0, \lim_{u\rw\IF}a_2(u)=-\IF$. Then by (\ref{y112})
$$y_i=x_i \mathbb{I} ( x_i>0, \lim_{u\rw\IF} a_i(u)=\IF\})-x_i \mathbb{I}( x_i>0, \lim_{u\rw\IF}a_i(u)=-\IF), i=1,2$$
we have
 $y_1=x_1$ and $y_2=-x_2$.
By \eqref{Bonfer}, in order to complete the proof, it suffices to prove that, as $u\rw\IF, S\rw\IF$, $\pi^-(u,S)\sim \pi^+(u,S)$ and to show that $\Sigma_{\re{i}}(u)=o(\pi^+(u,\re{S}))$, \re{$ i=1,2$}.\\
{\it \underline{Analysis of $\pi^{\pm}(u,S)$}}.
 We apply Lemma \ref{Pickands} to derive the uniform asymptotics for each term in the above sum. For this, we have to check conditions {\bf C0-C2}. Following the notation in Lemma \ref{Pickands}, let
$$\xi_{u,k}(t)=\overline{Z}_{u,k}(t), \quad \sigma_{\xi_{u,k}}=1, \quad g_k(u)=n^-_{u,k},\quad E=I_0, \quad  K_u=\{k: N_1(u)- 2\leq k\leq N_2(u)+2\}.$$
Conditions {\bf C0-C1} hold straightforwardly with $h(t)=0$. By (\ref{cor}), {\bf C2} holds {for $\zeta=\eta$, $\nu(u)=1$ and $\Delta(u)$ given in (\ref{delta(u)}). Hence we have that  $\zeta_\varphi(t)=\eta_\varphi(t)$}.
Thus by Lemma \ref{Pickands}, we have that for $S>x$
\BQN\label{uniform}
\lim_{u\rw\IF}\sup_{k\in K_u}\left|\frac{\pk{\int_{I_0} \mathbb{I}_0\left(\overline{Z}_{u,k}(t)-g_k(u)\right)\td t>x}}{\Psi(g_k(u))}-
\B_{\eta_\varphi}(x,[0,S])\right|=0.
\EQN
Consequently,
\BQN\label{upper}
\pi^{-}(u,S)&\re{=}& \sum_{k=N_1(u)- 2}^{N_2(u)+ 2}\pk{\int_{I_0} \mathbb{I}_0\left( \overline{Z}_{u,k}(t)-n^-_{u,k}\right)\td t>x}\nonumber\\
&\sim&\sum_{k=N_1(u)- 2}^{N_2(u)+ 2}\B_{\eta_\varphi}(x,[0,S])\Psi(n^-_{u,k})\nonumber\\
&\sim&
\B_{\eta_\varphi}(x,[0,S])\Psi(n(u))\sum_{k=N_1(u)- 2}^{N_2(u)+ 2}e^{-(1-\epsilon)n^2(u)\inf_{t\in I_k}w(g(u)|t|)}, \quad u\rw\IF.
\EQN
To get an upper bound \ehe{for} $\pi^{-}(u,S)$, it suffices to compute the sum above.  Note that $$\lim_{u\rw\IF}g(u)|a_i(u)|=0, \quad i=1,2.$$ Thus by Potter's theorem (see \cite{BI1989}) or Lemma 6.1 in \cite{KEP20151}, we have for any $0<\epsilon<\min(1,\beta)$ and all $u$  sufficiently large
\BQN\label{potter}
\frac{w(g(u)|s|)}{w(g(u)|t|)}\geq  (1-\epsilon/2)\min\left(\re{\left|\frac{s}{t}\right|^{\beta-\epsilon}, \left|\frac{s}{t}\right|^{\beta+\epsilon}} \right), \quad s,t\in E(u), t\neq 0
\EQN
implying that for $u$ large enough
\BQNY
\inf_{t\in I_k} w(g(u)|t|)\geq (1-\epsilon/2)\left(\frac{|k|-1}{|k|}\right)^{\beta-\epsilon}\sup_{t\in I_k}w(g(u)|t|), \quad  k\in K_u, k\neq 0.
\EQNY
Consequently, for any $0<\epsilon<1$, there exists $k_\epsilon\in \mathbb{N}$ such that for $|k|\geq k_\epsilon$ and $k\in K_u$  when $u$ is sufficiently large
\BQNY
\inf_{t\in I_k} w(g(u)|t|)\geq (1-\epsilon)\sup_{t\in I_k}w(g(u)|t|).
\EQNY
\K{Hence}, for $u$ large enough
\BQN\label{upper1}
\sum_{k=N_1(u)- 2}^{N_2(u)+ 2}e^{\re{-}(1-\epsilon)n^2(u)\inf_{t\in I_k}w(g(u)|t|)}&\leq& 2k_\epsilon+\sum_{k\in K_u, |k|\geq k_\epsilon }S^{-1}\int_{I_k}e^{-(1-\epsilon)^2n^2(u)w(g(u)|t|)}\td t\nonumber\\
&\leq& 2k_\epsilon+S^{-1}\int_{a_1(u)}^{a_2(u)}e^{-(1-\epsilon)^2n^2(u)w(g(u)|t|)}\td t\nonumber\\
&\leq& 2k_\epsilon+S^{-1}(g(u))^{-1}\int_{g(u)a_1(u)}^{g(u)a_2(u)}e^{-(1-\epsilon)^2n^2(u)w(|t|)}\td t.
\EQN
We have that for $S>x$
\BQN\label{upper2}
\int_{g(u)a_1(u)}^{g(u)a_2(u)}e^{-(1-\epsilon)^2n^2(u)w(|t|)}\td t\sim (1-\epsilon)^{-2/\beta} \inv{w}(n^{-2}(u))\beta^{-1}\int_{(1-\epsilon)^2y_1}^{(1-\epsilon)^2y_2} |t|^{1/\beta-1}e^{-|t|}\td t.
\EQN
The proof of (\ref{upper2}) is postponed {to} Appendix.
Further by (\ref{upper})-(\ref{upper2}) we obtain for $S>x$
\BQN\label{upper3}
\pi^{-}(u,S)\leq \frac{\B_{\eta_\varphi}(x,[0,S])}{S}\Theta(u)(1+o(1)), \quad u\rw\IF
\EQN
and similarly,
\BQN\label{lower}
\pi^{+}(u,S)\geq \frac{
\B_{\eta_\varphi}(x,[0,S])}{S}\Theta(u)(1+o(1)), \quad u\rw\IF.
\EQN
{\it \underline{Upper bound of $\Sigma_i(u), i=1,2$}}.
(\ref{uniform}) with $x=0$ gives that
\BQN\label{uniform1}
\lim_{u\rw\IF}\sup_{N_1(u)-2\leq k\leq N_2(u)+2}\left|\frac{\pk{\sup_{t\in I_k}\overline{Z}_{u}(t)>n^-_{u,k}}}{\Psi(n^-_{u,k})}-\B_{\eta_\varphi}(0,[0,S])\right|=0.
\EQN
Thus in light of (\ref{upper})-(\ref{upper3}) we have that
\BQN\label{neg}
\Sigma_1(u)&\leq&\sum_{N_1(u)-2\leq k\leq N_2(u)+1}\left(\pk{\sup_{t\in I_k}\overline{Z}_{u}(t)>n^-_{u,k}} +\pk{\sup_{t\in I_{k+1}}\overline{Z}_{u}(t)>n^-_{u,k+1}}\right.\nonumber\\
&&\quad \left.-\pk{\sup_{t\in I_k\cup I_{k+1}}\overline{Z}_{u}(t)>\widetilde{n}_{u,k}}\right)\nonumber\\
&\leq & \sum_{N_1(u)-2\leq k\leq N_2(u)+1}\B_{\eta_\varphi}(0,[0,S])\Psi(n^-_{u,k})(1+o(1))\\
&&  +\sum_{N_1(u)-2\leq k\leq N_2(u)+1} \B_{\eta_\varphi}(0,[0,S]) \Psi(n^-_{u,k+1})(1+o(1))\nonumber\\
&&\quad -\sum_{N_1(u)-2\leq k\leq N_2(u)+1}\B_{\eta_\varphi}(0,[0,2S])\Psi(\widetilde{n}_{u,k})(1+o(1))\nonumber\\
&\leq& \mathbb{Q}S^{-1}\left(2 \B_{\eta_\varphi}(0,[0,S])-\B_{\eta_\varphi}(0,[0,2S])\right)\Theta(u), \quad u\rw\IF,
\EQN
where $\widetilde{n}_{u,k}=\max(n^-_{u,k}, n^-_{u,k+1})$.
Moreover, by \cite{KEP2016}[Corollary 3.1] and (\ref{cor}) there exists $\mathcal{C}, \mathcal{C}_1>0$ such that for
$N_1(u)-2\leq k<k+1<l\leq N_2(u)+1$ and all $u$ large enough
\BQNY
\pk{\sup_{t\in I_k}\overline{Z}_{u}(t)>n^-_{u,k}, \sup_{t\in I_{l}}\overline{Z}_{u}(t)>n^-_{u,l}}\leq \mathcal{C}S^2e^{-\mathcal{C}_1|k-l|^\gamma S^\gamma}\Psi(\re{\hat{n}_{u,k,l}}),
\EQNY
with $\hat{n}_{u,k,l}=\min(n^-_{u,k}, n^-_{u,l})$ and $\gamma=\min(\alpha_0, \alpha_\IF)$,
which combined with (\ref{upper})-(\ref{upper3}) leads to
\BQN\label{neg1}
\Sigma_2(u)&\leq& \sum_{N_1(u)-2\leq k<k+1<l\leq N_2(u)+1}\mathcal{C}S^2e^{-\mathcal{C}_1|k-l|^\gamma S^\gamma}\Psi(\hat{n}_{u,k,l})\nonumber\\
&\leq& \sum_{N_1(u)-2\leq k\leq N_2(u)+1}\Psi(n^-_{u,k})\sum_{l\geq \re{2}}\mathcal{C}S^2e^{-\mathcal{C}_1l^\gamma S^\gamma}\nonumber\\
&\leq& \sum_{N_1(u)-2\leq k\leq N_2(u)+1}\Psi(n^-_{u,k})\mathbb{Q}S^2e^{-\mathbb{Q}_1 S^\gamma}\nonumber\\
&\leq& \mathbb{Q}S^2e^{-\mathbb{Q}_1 S^\gamma}\Theta(u), \quad u\rw\IF.
\EQN
Consequently, by (\ref{neg}) and (\ref{neg1}), for any $S>0$
\BQN\label{neg2}
\Sigma_1(u)+\Sigma_2(u)\leq \mathbb{Q}\left(S^{-1}\left(2\B_{\eta_\varphi}(0,[0,S])-\B_{\eta_\varphi}(0,[0,2S])\right)+S^2e^{-\mathbb{Q}_1 S^\gamma}\right)\Theta(u).
\EQN
{\it \underline{Exact asymptotics of $\pi(u)$}}.
Inserting (\ref{upper3}), (\ref{lower}) and (\ref{neg2}) into (\ref{Bonfer}) and dividing each term by $\Theta(u)$, we have that for $S>x$
\BQN
 \limsup_{u\rw\IF}\frac{\pi(u)}{\Theta(u)}&\leq& \frac{\B_{\eta_\varphi}(x,[0,S])}{S}+\mathbb{Q}\left(2\frac{\B_{\eta_\varphi}(0,[0,S])}{S}-\frac{\B_{\eta_\varphi}(0,[0,2S])}{S}+S^2e^{-\mathbb{Q}_1 S^\gamma}\right)\label{final}\\
 \liminf_{u\rw\IF}\frac{\pi(u)}{\Theta(u)}&\geq&  \frac{\B_{\eta_\varphi}(x,[0,S])}{S}-\mathbb{Q}\left(2\frac{\B_{\eta_\varphi}(0,[0,S])}{S}-\frac{\B_{\eta_\varphi}(0,[0,2S])}{S}+S^2e^{-\mathbb{Q}_1 S^\gamma}\right).\label{final1}
\EQN
Combination of (\ref{final})-(\ref{final1}) and the fact that (see \cite{Pit96} and \cite{DI2005})
$$\lim_{S\rw\IF}\frac{\B_{\eta_\varphi}(0,[0,S])}{S}\in (0,\IF)$$
leads to
$$\liminf_{S\rw\IF}\frac{\B_{\eta_\varphi}(\ehe{x},[0,S])}{S}=\limsup_{S\rw\IF}\frac{\B_{\eta_\varphi}(\ehe{x},[0,S])}{S}<\IF.$$
We next prove that
\bqn{ \label{ns3}
\liminf_{S\rw\IF}\frac{\B_{\eta_\varphi}(\ehe{x},[0,S])}{S}>0.
}
 Note that for $\pk{\int_{E(u)\bigcap (\cup_{k\in \mathbb{Z}} I_{2k})} \mathbb{I}_0(Z_u(t)-n(u))\td t>x}$, similarly as in (\ref{final})-(\ref{final1}) we have for $S>x$
 \re{
\BQNY
 \liminf_{u\rw\IF}\frac{\pk{\int_{E(u) \bigcap (\cup_{k\in \mathbb{Z}} I_{2k})} \mathbb{I}_0(Z_u(t)-n(u))\td t>x}}{\Theta(u)}&\geq& \frac{\B_{\eta_\varphi}(x,[0,S])}{2S}-\mathbb{Q}S^2e^{-\mathbb{Q}_1 S^\gamma},
\EQNY which together with (\ref{final}) gives that   that for any $ S_1>x$
\BQNY
\liminf_{S\rw\IF}\frac{\B_{\eta_\varphi}(\ehe{x},[0,S])}{S}\geq \frac{\B_{\eta_\varphi}(x,[0,S_1])}{S_1}-\mathbb{Q}S_1^2e^{-\mathbb{Q}_1 S_1^\gamma}.
\EQNY}
Notice  that for $S>x$,
 $\B_{\eta_\varphi}([0,S],x)$ is non-decreasing with respect to $S$ for $S>x$ and
   $$\B_{\eta_\varphi}([0,S],x)\geq \int_{\mathbb{R}}\pk{\inf_{t\in [0,S]}\left( \sqrt{2}\eta_\varphi(t)-Var(\eta_\varphi(t))\right)>-z}e^{-z}\td z>0.$$
   Hence for some $S_1>x$,
\BQNY
\liminf_{S\rw\IF}\frac{\B_{\eta_\varphi}(x,[0,S])}{S}\geq \frac{\B_{\eta_\varphi}(x,[0,S_1])}{S_1}-\mathbb{Q}S_1^2e^{-\mathbb{Q}_1 S_1^\gamma}>0,
\EQNY
which confirms that (\ref{ns3}) holds. Hence we have
$$\B_{\eta_\varphi}(x) =\lim_{S\rw\IF}\frac{\B_{\eta_\varphi}(x,[0,S])}{S}\in (0,\IF).$$
Letting $S\rw\IF$ in (\ref{final})-(\ref{final1}), we derive
$$\pi(u)\sim \B_{\eta_\varphi}(x) \Theta(u), \quad u\rw\IF.$$
This completes the proof \kk{of case i)}.\\
\underline{$\diamond$ \it Case ii)} Let us first  assume that
$$\lim_{u\rw\IF}a_1(u)=-\IF, \quad \lim_{u\rw\IF}a_2(u)=\IF.$$
 By (\ref{var}) observe that  for $u$ sufficiently large
\BQN\label{Bonfer1}
p_0(u)\leq \pi(u)\leq p_0(u)+\pi_1(u),
\EQN
where
\BQNY
p_0(u)=\pk{\int_{[-S,S]} \mathbb{I}_0(Z_u(t)-n(u))\td t>x},\quad
\pi_1(u)=\sum_{k=N_1(u)- 2,k\neq -1, 0 }^{N_2(u)+2}\pk{\sup_{{t\in I_0}}\overline{Z}_{u,k}(t)>n_{u,k}^{-}}.
\EQNY
In order to complete the proof, we shall derive the asymptotics of $p_0(u)$ and show further that  $$\pi_1(u)=o(p_0(u)), \quad u \rw\IF, S\rw\IF.$$
{\it \underline{The analysis of $p_0(u)$}}.
In order to apply Lemma \ref{Pickands},   we need to check the validity of {\bf C0-C2}. {\bf C0} holds straightforwardly. By (\ref{gammalim}) and uniform convergence theorem (e.g., in \cite{BI1989})
\BQNY
&&\lim_{u\rw\IF}\sup_{t\in[-S,S]} \left|n^2(u)w(g(u)|t|)- \gamma|t|^\beta\right|\\
&&\leq \lim_{u\rw\IF}\sup_{t\in[-S,S]} \left|\left(n^2(u)w(g(u))-\gamma\right)\frac{w(g(u)|t|)}{w(g(u))}\right| +\gamma\lim_{u\rw\IF}\sup_{t\in[-S,S]}\left|\frac{w(g(u)|t|)}{w(g(u))}- |t|^\beta\right|=0.
\EQNY
 By (\ref{var}), (\ref{gammalim}) and uniform convergence theorem (e.g., in \cite{BI1989}), we have that for any $S>0$
 \BQNY
 \lim_{u\rw\IF}\sup_{t\in[-S,S]}\left|n^2(u)(1-\sigma_u(t))-\gamma |t|^\beta\right|&\leq& \lim_{u\rw\IF}\sup_{t\in[-S,S]}  \left|n^2(u)w(g(u)|t|)\right|\left| \frac{1-\sigma_u(t)}{w(g(u)|t|)}-1\right|\\
 && \quad +\lim_{u\rw\IF}\sup_{t\in[-S,S]} \left|n^2(u)w(g(u)|t|)- \gamma|t|^\beta\right|\\
 &=&0,
 \EQNY
 which confirms that {\bf C1} hold with $h(t)=\gamma |t|^\beta$.
By (\ref{cor}), {\bf C2} is satisfied with  \re{ $\zeta_\varphi(t)=\eta_\varphi(t)$}.
 Thus  we have
 \BQN\label{p0}
 p_0(u)\sim \B_{\eta_\varphi}^{\gamma|t|^\beta}(x,[-S,S])\Psi(n(u)).
 \EQN
 {\it \underline{Upper bound of $\pi_1(u)$}}.
\ehe{By (\ref{uniform})} 
\BQNY
\pi_1(u)&\sim& \sum_{k=N_1(u)- 2,k\neq -1, 0 }^{N_2(u)+2} \B_{\eta_\varphi}(0,[0,S])\Psi(n_{u,k}^{-})\\
&\leq & \B_{\eta_\varphi}(0,[0,S])\Psi(n(u))\sum_{k=N_1(u)- 2,k\neq -1, 0 }^{N_2(u)+2} e^{-(1-\epsilon)\inf_{t\in I_k} n^2(u)w(g(u)|t|)}, \quad u\rw\IF.
\EQNY
Further, using  (\ref{potter}) for  $u$ sufficiently large  and $S>1$ we have
\BQNY
\inf_{t\in I_k} n^2(u)w(g(u)|t|)\geq \frac{\gamma}{2}\inf_{t\in I_k} \frac{w(g(u)|t|)}{w(g(u))}\geq \mathbb{Q}(|k|S)^{\beta/2}, \quad  N_1(u)- 2\leq k\leq N_2(u)+2, k\neq -1, 0.
\EQNY
Thus we have
\BQN\label{pi1}
\pi_1(u)
&\leq & \B_{\eta_\varphi}(0,[0,S])\Psi(n(u))\sum_{k=N_1(u)- 2,k\neq -1, 0 }^{N_2(u)+2} e^{- \mathbb{Q}(|k|S)^{\beta/2}}\nonumber\\
&\leq &   \B_{\eta_\varphi}(0,[0,S]) \re{\mathbb{Q}_2 e^{-\mathbb{Q}_1S^{\beta/2}}}\Psi(n(u)), \quad u\rw\IF.
\EQN
{\it \underline{Exact asymptotics of $\pi(u)$}}.
Inserting (\ref{p0}) and (\ref{pi1}) into (\ref{Bonfer1}) and dividing each terms by $\Psi(n(u))$, we have
\BQN\label{uplo}
\B_{\eta_\varphi}^{\gamma|t|^\beta}(x,[-S,S])&\leq& \re{\liminf}_{u\rw\IF}\frac{\pi(u)}{\Psi(n(u))}\nonumber\\
&\leq& \re{\limsup}_{u\rw\IF} \frac{\pi(u)}{\Psi(n(u))}\nonumber\\
&\leq& \B_{\eta_\varphi}^{\gamma|t|^\beta}(x, [-S,S])+\B_{\eta_\varphi}(0,[0,S]) \re{\mathbb{Q}_2}e^{-\mathbb{Q}_{\re{1}}S^{\beta/2}}.
\EQN
This implies that for $S_1>x$
\BQNY
0<\lim_{S\rw\IF}\B_{\eta_\varphi}^{\gamma|t|^\beta}(\ehe{x},[-S,S])\leq \B_{\eta_\varphi}^{\gamma|t|^\beta}(\ehe{x},[-S_1,S_1])+\B_{\eta_\varphi}(0,[0,S_1])\re{\mathbb{Q}_2}e^{-\mathbb{Q}_{\re{1}}S_1^{\beta/2}}<\IF.
\EQNY
Moreover, by case i)
$$\lim_{S\rw\IF}\frac{\B_{\eta_\varphi}(0,[0,S])}{S}\in (0,\IF).$$
Therefore, letting $S\rw \IF$ in (\ref{uplo}), we have
$$\pi(u)\sim \B_{\eta_\varphi}^{\gamma|t|^\beta}(x,(-\IF,\IF))\Psi(n(u)), \quad u\rw\IF.$$
This establishes the claim for $a_1=-\IF, a_2=\IF$.\\
Next we focus on the case $a_i\in (-\IF, \IF), i=1,2$ with $a_2-a_1>x$. Note that for $\epsilon>0$ sufficiently small and $u$ sufficiently large
\BQNY
\pk{\int_{[a_1+\epsilon,a_2-\epsilon]} \mathbb{I}_0(Z_u(t)-n(u))\td t>x}\leq \pi(u)\leq \pk{\int_{[a_1-\epsilon,a_2+\epsilon]} \mathbb{I}_0(Z_u(t)-n(u))\td t>x}.
\EQNY
Applying Lemma \ref{Pickands}
it follows that
\BQN\label{caseii}
 \B_{\eta_\varphi}^{\gamma|t|^\beta}(x,[a_1+\epsilon,a_2-\epsilon])\Psi(n(u))\leq \pi(u)\leq  \B_{\eta_\varphi}^{\gamma|t|^\beta}(x,[a_1-\epsilon,a_2+\epsilon])\Psi(n(u)), \quad u\rw\IF.
 \EQN
Notice that
\BQNY\B_{\eta_\varphi}^{\gamma|t|^\beta}(x,[a_1-\epsilon,a_2+\epsilon])&=&\int_{\mathbb{R}}\pk{\int_{[a_1-\epsilon,a_2+\epsilon]} \mathbb{I}_0(\widetilde{\eta}(t)+z)\td t>x}e^{-z}\td z\\
&=&\int_{\mathbb{R}}\pk{\int_{[a_1,a_2]} \mathbb{I}_0(\widetilde{\eta}(t)+z)\td t+\int_{[a_1-\epsilon,a_1]\cup[a_2, a_2+\epsilon] } \mathbb{I}_0(\widetilde{\eta}(t)+z)\td t >x}e^{-z}\td z,
\EQNY
where $\widetilde{\eta}(t)=\sqrt{2}\eta_\varphi(t)-Var(\eta_\varphi(t))-\gamma|t|^\beta$, and
$$0\leq \int_{[a_1-\epsilon,a_1]\cup[a_2, a_2+\epsilon] } \mathbb{I}_0(\widetilde{\eta}(t)+z)\td t\leq 2\epsilon.$$
Hence
\BQNY\B_{\eta_\varphi}^{\gamma|t|^\beta}(x,[a_1,a_2])\leq \B_{\eta_\varphi}^{\gamma|t|^\beta}(x,[a_1-\epsilon,a_2+\epsilon])\leq \B_{\eta_\varphi}^{\gamma|t|^\beta}(x-2\epsilon,[a_1,a_2]).
\EQNY
The continuity of $\B_{\eta_\varphi}^{\gamma|t|^\beta}(x,[a_1,a_2])$ ( see  Lemma \ref{Pickands}) leads to
$$\lim_{\epsilon\rw 0} \B_{\eta_\varphi}^{\gamma|t|^\beta}(x-2\epsilon,[a_1,a_2])=\B_{\eta_\varphi}^{\gamma|t|^\beta}(x,[a_1,a_2]),$$
implying that
\BQN\label{continue1}\lim_{\epsilon\rw 0} \B_{\eta_\varphi}^{\gamma|t|^\beta}(x,[a_1-\epsilon,a_2+\epsilon])=\B_{\eta_\varphi}^{\gamma|t|^\beta}(x,[a_1,a_2]).
\EQN
We can analogously show that
\BQN\label{continue2}\lim_{\epsilon\rw 0} \B_{\eta_\varphi}^{\gamma|t|^\beta}(x,[a_1+\epsilon,a_2-\epsilon])=\B_{\eta_\varphi}^{\gamma|t|^\beta}(x,[a_1,a_2]).
\EQN
Letting $\epsilon\rw 0$ in (\ref{caseii}), we have that
$$\pi(u)\sim  \B_{\eta_\varphi}^{\gamma|t|^\beta}(x,[a_1,a_2])\Psi(n(u)), \quad u\rw\IF.$$
This establishes the claim for case $a_i\in (-\IF, \IF), i=1,2$ with $a_2-a_1>x$.\\
  For the case that $a_1=-\IF$, $a_2\in (-\IF,\IF)$, and $a_1\in (-\IF,
  \IF)$, $a_2=\IF$,  we can establish the claim by using same approach.\\
{\underline{$\diamond$ \it Case iii)}} Let us first consider the case that $$\lim_{u\rw\IF}a_1(u)\frac{g(u)}{\inv{w}(n^{-2}(u))}=-\IF, \quad  \lim_{u\rw\IF}a_2(u)\frac{g(u)}{\inv{w}(n^{-2}(u))}=\IF.$$
 By (\ref{var}) we have
\BQN\label{Bonfer2}
p_1(u)\leq \pi_2(u)\leq \sum_{i=1}^3p_i(u)+\pi_3(u),
\EQN
where
\BQNY
\pi_2(u)&=&\pk{\int_{E(u)} \mathbb{I}_0(Z_u(t)-n(u))\td t>\theta(u)x}\\
p_1(u)&=&\pk{\int_{[-\theta(u) S,\theta(u) S]} \mathbb{I}_0(Z_u(t)-n(u))\td t>\theta(u)x},\\
p_i(u)&=& \pk{\sup_{t\in I_{i-3}\setminus [-\theta(u) S,\theta(u) S]}\overline{Z}_{u}(t)>n'(u)}, \quad i=2,3,\\
\pi_3(u)&=&\sum_{k=N_1(u)-2, k\neq 0, -1}^{N_2(u)+2} \pk{\sup_{t\in I_k} \overline{Z}_{u}(t)>n_{u,k}^-},
\EQNY
with
$${\theta(u)}=\frac{\inv{w}(n^{-2}(u))}{g(u)},\quad   \re{n'(u)}=n(u)\left(1+(1-\epsilon)\inf_{t\in [-S,S]\setminus[-\theta(u) S, \theta(u) S]}w(g(u)|t|)\right), \quad 0<\epsilon<1.$$
We shall derive the exact asymptotics of $p_1(u)$ and then prove that $p_2(u)$, $p_3(u)$ and $\pi_3(u)$ are all negligible compared with $p_1(u)$ as $u\rw\IF$ and $S\rw\IF$.\\
{\it \underline{Analysis of $p_1(u)$}}.
Substituting $t$ by $\theta(u)s$ we obtain
\BQNY
p_1(u) =\pk{\int_{[-S,S]} \mathbb{I}_0(Z_u(\theta(u)t)-n(u))\td t>x}.
\EQNY
Next we check {\bf C0}-{\bf C2} in Lemma \ref{Pickands}. {\bf C0} holds straightforwardly. Moreover,  by (\ref{var}) and uniform convergence theorem (e.g., in \cite{BI1989}), with noting that $n^{2}(u)=(w(\inv{w}(n^{-2}(u))))^{-1}$, we have that for any $S>0$
 \BQNY
 \lim_{u\rw\IF}\sup_{t\in[-S,S]}\left|n^2(u)(1-\sigma_u(\theta(u)t))- |t|^\beta\right|&\leq& \lim_{u\rw\IF}\sup_{t\in[-S,S]}  \left|\frac{w(\inv{w}(n^{-2}(u))|t|)}{w(\inv{w}(n^{-2}(u)))}\right|\left| \frac{1-\sigma_u(\theta(u)t)}{w(\inv{w}(n^{-2}(u))|t|)}-1\right|\\
 && \quad +\lim_{u\rw\IF}\sup_{t\in[-S,S]} \left|\frac{w(\inv{w}(n^{-2}(u))|t|)}{w(\inv{w}(n^{-2}(u)))}- |t|^\beta\right|\\
 &=&0.
 \EQNY
 This confirms that {\bf C1} holds with $h(t)=|t|^\beta$. It follows from (\ref{cor}) that
$$\lim_{u\rw\IF}\sup_{s,t\in [-S,S], s\neq t}\left|\frac{n^2(u)(1-Corr(\re{\overline{Z}_u}(\theta(u)t)-\re{\overline{Z}_u}(\theta(u)s)))}{ \frac{\sigma_\eta^2(\theta(u)\Delta(u)|t-s|)}
{\sigma_\eta^2(\Delta(u))}}-1\right|=0$$
with
$$\lim_{u\rw\IF}\theta(u)=\lim_{u\rw\IF}\frac{\inv{w}(n^{-2}(u))}{g(u)}=0.$$
This means that {\bf C2} is satisfied with $\zeta_\varphi=0$.
Therefore we have that
\BQN\label{p1}
p_1(u)\sim \B_0^{|t|^\beta}(x,[-S,S])\Psi(n(u)).
\EQN
{\it \underline{Upper bound for $\pi_3(u)$, $p_i(u), i=2,3$}}.
Next we find an upper bound of $\pi_3(u)$. Similarly as (\ref{pi1}), by (\ref{uniform}), we have
\BQN\label{pi2}
\pi_3(u)\leq \B_{\eta_\varphi}(0,[0,S])   \mathbb{Q}_2e^{-\mathbb{Q}_{1}S^{\beta/2}}\Psi(n(u)), \quad u\rw\IF.
\EQN
Finally, we focus on deriving an upper bound of  $p_i(u), i=2,3$. In light of (\ref{uniform}), we have that
\BQN\label{pi}
p_i(u)\leq  \pk{\sup_{t\in I_{i-3}}\re{\overline{Z}}_{u}(t)>n'(u)}\sim \B_{\eta_\varphi}(0,[0,S])\Psi(\re{n'(u)}), \quad i=2,3.
\EQN
By (\ref{potter}), we have that
\BQNY
n^2(u)\inf_{t\in [-S,S]\setminus[-\theta(u) S,\theta(u)) S]}w(g(u)|t|)\geq \re{\frac{1}{2}} \inf_{t\in [0, S/\theta(u)]\setminus[- S, S]}\frac{w(\inv{w}(n^{-2}(u))|t|)}{w(\inv{w}(n^{-2}(u)))}\geq \frac{1}{4} \inf_{t\in [0, S/\theta(u)]\setminus[- S, S]}|t|^{\beta/2}\geq \frac{1}{4}|S|^{\beta/2},
\EQNY
which together with (\ref{pi}) leads to
\BQN\label{piu}
p_i(u)\leq \B_{\eta_\varphi}(0,[0,S]) e^{-\mathbb{Q}|S|^{\beta/2}/8}\Psi(n(u)), \quad  u\rw\IF, i=2,3.
\EQN
{\it \underline{Exact asymptotics of $\pi_2(u)$}}.
Inserting (\ref{p1}), (\ref{pi2}) and  (\ref{piu}) into (\ref{Bonfer2}), we have that
\BQN\label{uplo1}
\liminf_{u\rw\IF}\frac{\pi_2(u)}{\Psi(n(u))}&\geq & \mathcal{B}_0^{|t|^\beta}(x,[-S,S]) , \nonumber \\
\limsup_{u\rw\IF}\frac{\pi_2(u)}{\Psi(n(u))}&\leq& \mathcal{B}_0^{|t|^\beta}(x,[-S,S])+\B_{\eta_\varphi}(0,[0,S]) e^{-\mathbb{Q}S^{\beta/2}},
\EQN
which implies
that
\BQNY
0<\lim_{S\rw\IF}\mathcal{B}_0^{|t|^\beta}(x,[-S,S])\leq \mathcal{B}_0^{|t|^\beta}(x,[-S_1,S_1])+\B_{\eta_\varphi}(0,[0,S_1])e^{-\mathbb{Q}S_1^{\beta/2}}<\IF, \quad S_1>x.
\EQNY
Letting $S\rw\IF$ in (\ref{uplo1}) yields that
$$\pi_2(u)\sim \mathcal{B}_0^{|t|^\beta}(x,(-\IF,\IF))\Psi(n(u)).$$
This establishes the claim for $b_1=-\IF$, $b_2=\IF$.\\ Next we focus on the case $b_1, b_2\in (-\IF,\IF)$ with $b_2-b_1>x$.  For $\epsilon>0$ sufficiently small and $u$ sufficiently large, we have
$$\pk{\int_{[\theta(u) (b_1+\epsilon),\theta(u) (b_2-\epsilon)]} \mathbb{I}_0(Z_u(t)-n(u))\td t>\theta(u)x}\leq \pi_2(u)\leq \pk{\int_{[\theta(u) (b_1-\epsilon),\theta(u) (b_2+\epsilon)]} \mathbb{I}_0(Z_u(t)-n(u))\td t>\theta(u)x}.$$
Applying Lemma \ref{Pickands}, similarly as in (\ref{p1}), we derive that
$$\mathcal{B}_0^{|t|^\beta}(x,[b_1+\epsilon,b_2-\epsilon])\Psi(n(u))\leq \pi_2(u)
\leq \mathcal{B}_0^{|t|^\beta}(x,[b_1-\epsilon,b_2+\epsilon])\Psi(n(u)), \quad u\rw\IF.$$
In light of (\ref{continue1}) and (\ref{continue2}), letting $\epsilon\rw 0$ in above inequality, we get
$$\pi_2(u)
\sim \mathcal{B}_0^{|t|^\beta}(x,[b_1,b_2])\Psi(n(u)), \quad u\rw\IF.$$
Hence the claim for case $b_1, b_2\in (-\IF,\IF)$ is established.\\
{The proofs of claims for other cases of $b_1$ and $b_2$ can be done by the same line of reasoning.}

This completes the proof.\QED
\prooftheo{TH2} We set $\Delta(u)=\overleftarrow{\sigma}
 \left(\frac{\sqrt{2}\sigma^2(ut^*)}{u(1+ct^*)}\right)$. By (\ref{scale}) observe that for any $u>0$ (recall the definition of
$\overleftarrow{\sigma}$ in \eqref{zhr})
\BQN\label{Bonfer3}
\pi_4(u)\leq\pk{\frac{1}{\Delta(u)}\int_0^\IF \mathbb{I}_u(X(t)-ct)\td t>x}\leq
\pi_4(u)+\pi_5(u)
\EQN
where
\BQN\label{pi4}\pi_4(u)&=&\pk{\frac{u}{\Delta(u)}\int_{E_1(u)} \mathbb{I}_0\left(\frac{X(ut)}{u(1+ct)}M(u)-M(u)\right)\td t>x},\nonumber\\
 \pi_5(u)&=&\pk{ \sup_{t\in [0,\IF]\setminus E_1(u)}\frac{X(ut)}{u(1+ct)}M(u)>M(u)},
\EQN
with
\BQN\label{yu}
\quad E_1(u)=\left[t_u-\frac{\ln M(u)}{M(u)}, t_u+\frac{\ln M(u)}{M(u)}\right].
\EQN
The idea of the proof is \re{to derive the exact asymptotics of $\pi_4(u)$ using Theorem \ref{TH1} } and to show that $\pi_5(u)=o(\pi_4(u))$ as $u\rw\IF$.\\
{\it \underline{ The analysis of $\pi_4(u)$}}. Scaling time with $\frac{u}{\Delta(u)}$, we have
\BQN\label{pi3}\pi_4(u)=\pk{\int_{E_2(u)} \mathbb{I}_0\left(Z_u(t)-M(u)\right)\td t>x},
\EQN
where
 \BQN\label{zu}
Z_u(t)=\frac{X(ut_u+\Delta(u)t)}{u(1+ct_u)+c\Delta(u)t}, \quad E_2(u)=\left[-\frac{u\ln M(u)}{\Delta(u)M(u)}, \frac{u\ln M(u)}{\Delta(u)M(u)}\right].
 \EQN
It follows from  \eqref{L1} and \eqref{L2} that
\BQN\label{corzu}
&&1-\sqrt{Var(Z_u(t))}\sim \left(\sqrt{\frac{B}{2A}}\frac{\Delta(u)}{u}|t|\right)^2, \quad t\in E_2(u),\nonumber \\
 &&\lim_{u\rw \IF}\sup_{s\neq t, s,t \in E_2(u)}\left|\re{M^2(u)}\frac{1-Corr(Z_u(s), Z_u(t))}{\frac{\sigma^2(\Delta(u)|s-t|)}{\re{\sigma^2(\Delta(u))}}}-1\right|=0,
\EQN
which imply that  (\ref{var}) and (\ref{cor}) hold with
\BQN\label{notation}
w(t)&=&t^2, \quad g(u)=\sqrt{\frac{B}{2A}}\frac{\Delta(u)}{u}, \quad n(u)=M(u), \quad  \eta=X, \quad \Delta(u)=\overleftarrow{\sigma}
 \left(\frac{\sqrt{2}\sigma^2(ut^*)}{u(1+ct^*)}\right), \nonumber\\ a_1(u)&=&-\frac{u\ln M(u)}{\Delta(u)M(u)},\quad
a_2(u)=\frac{u\ln M(u)}{\Delta(u)M(u)}.
\EQN
Next we check the assumptions of i) in Theorem \ref{TH1}.  Following the notation in Theorem \ref{TH1}, we have
\BQN\label{condition11}
\lim_{u\rw\IF}g(u)=\lim_{u\rw\IF}\sqrt{\frac{B}{2A}}\frac{\Delta(u)}{u}=0,  \quad \lim_{u\rw\IF}g(u)|a_i(u)|=\lim_{u\rw\IF}\sqrt{\frac{B}{2A}}\frac{\ln M(u)}{M(u)}=0, i=1,2.
\EQN
Note that
$$n^2(u)w(g(u))=\frac{B}{2A}\left(\frac{M(u)\Delta(u)}{u}\right)^2\sim \mathbb{Q}\left(\frac{\inv{\sigma}(u^{-1}\sigma^2(u))}{\sigma(u)}\right)^2 $$
is a regularly varying function at $\IF$ with index $2\tau$
with
\BQN\label{re}
\tau=\left\{\begin{array}{cc}
\frac{2\alpha_\IF-1}{\alpha_0}-\alpha_\IF, & \varphi=0\\
-\alpha_\IF, & \varphi\in (0,\IF)\\
\frac{2\alpha_\IF-1}{\alpha_\IF}-\alpha_\IF, & \varphi=\IF.
\end{array}\right.
\EQN
Since $\tau<0$ for all $\varphi$, then for all $\varphi\in [0,\IF]$
 \BQN\label{condition12}
\gamma=\lim_{u\rw\IF}n^2(u)w(g(u))=0.
\EQN
 Moreover,
\BQN\label{y11}\lim_{u\rw\IF}n^2(u)w(g(u)|a_{\re{i}}(u)|)=\lim_{u\rw\IF} \frac{B}{2A}(\ln M(u))^2=\IF, \quad i=1, 2
\EQN
  and by (\ref{re})
\BQN\label{y12}\lim_{u\rw\IF} a_1(u)=-\lim_{u\rw\IF}\frac{u\ln M(u)}{\Delta(u)M(u)}=-\mathbb{Q}\lim_{u\rw\IF}\frac{\sigma(u)\re{\ln M(u)}}{\inv{\sigma}(u^{-1}\sigma^2(u))}=-\IF,\EQN
   which imply that $y_1=-\IF$. Similarly we can check that $y_2=\IF$. Additionally,
   \BQN\label{condition13}
   \lim_{u\rw\IF}n(u)w(g(u)|a_i(u)|)=\frac{B}{A}\lim_{u\rw\IF}\frac{(\ln M(u))^2}{M(u)}=0, \quad i=1,2.
   \EQN
Hence all the assumptions  in i) of Theorem \ref{TH1} are satisfied, leading to
\BQN\label{upper4}
\pi_4(u)&\sim& \B_{X_{\varphi}}(x)\frac{1}{2}\int_{-\IF}^{\IF}|t|^{1/2-1}e^{-|t|}\td t\sqrt{\frac{2A}{B}}\frac{u}{M(u)\Delta(u)}\Psi(M(u))\nonumber\\
&\sim& \B_{X_{\varphi}}(x)\sqrt{\frac{2A\pi}{B}} \frac{u}{M(u)\Delta(u)}\Psi(M(u)), \quad u\rw\IF,
\EQN
{where $X_{\varphi}$ is given in (\ref{theta}).}

{\it \underline{Upper bound of $\pi_5(u)$}}.
By  \cite{DI2005}[Lemma 7] or  \cite{KrzysPeng2015}[Lemma 5.6], we have
\BQN\label{pi4upper}
\pi_5(u)=o\left(\frac{u}{M(u)\Delta(u)}\Psi(M(u))\right), \quad u\rw\IF,
\EQN
which combined with (\ref{Bonfer3}) and (\ref{upper4}) leads to
\BQNY
\pk{\frac{1}{\Delta(u)}\int_0^\IF \mathbb{I}_u(X(t)-ct)\td t>x}\sim\B_{X_{\varphi}}(x)\sqrt{\frac{2A\pi}{B}} \frac{u}{M(u)\Delta(u)}\Psi(M(u)), \quad u\rw\IF.
\EQNY
This completes the proof.\QED
\proofkorr{TH3} We also set $ \Delta(u)=\overleftarrow{\sigma}
 \left(\frac{\sqrt{2}\sigma^2(ut^*)}{u(1+ct^*)}\right)$. Observe that for $0\leq x_1\leq x_2<\IF$ and $u>0$ the conditional distribution can be rewritten  as the ratio of two sojourn probabilities  over different intervals, i.e.,
\BQN\label{ratio}
\pk{\frac{\tau_u^*(x_1,x_2)-ut_u}{A(u)}<y}=\frac{q(u)}{\pk{\frac{1}{\Delta(u)}\int_0^\IF \mathbb{I}_u(X(t)-ct)\td t>x_1}},
\EQN
where {$y\in\mathbb{R}$ and}
$$q(u)=\pk{\frac{1}{\Delta(u)}\int_{[0, ut_u+A(u)y]} \mathbb{I}_u(X(t)-ct)\td t>x_2}.$$
Hence in order to complete the proof it suffices to derive the asymptotics of $q(u)$. Using  notation for $ E_1(u)$ and $\pi_5(u)$ introduced in (\ref{pi4}) and (\ref{yu}), we have that
  \BQN\label{Bonfer4}
  q_1(u)\leq q(u)\leq q_1(u)+\pi_4(u),
  \EQN
  with
  $$q_1(u)=\pk{\frac{u}{\Delta(u)}\int_{\left[t_u-\frac{\ln M(u)}{M(u)}, t_u+\frac{A(u)}{u}y\right]} \mathbb{I}_0\left(\frac{X(ut)}{u(1+ct)}M(u)-M(u)\right)\td t>x_2}.$$
  Scaling  time by $\frac{u}{\Delta(u)}$, we rewrite $$q_1(u)=\pk{\int_{E_3(u)} \mathbb{I}_0\left(Z_u(t)-M(u)\right)\td t>x_2},$$
where
$$ Z_u(t)=\frac{X(ut_u+\Delta(u)t)}{u(1+ct_u)+c\Delta(u)t}, \quad E_3(u)=\left[-\frac{u\ln M(u)}{\Delta(u)M(u)}, \frac{A(u)}{\Delta(u)}y\right].$$
 Using the same notation as introduced in (\ref{notation}) with the exception that  $a_2(u)=\frac{A(u)}{\Delta(u)}y$, we have that (\ref{corzu}), (\ref{condition11}), (\ref{condition12}) and (\ref{condition13}) also hold. We next get the value of $y_i, i=1,2$. For this, note that
 \BQNY
\lim_{u\rw\IF}n^2(u)w(g(u)|a_2(u)|)=\lim_{u\rw\IF} \frac{B}{2A}\left(\frac{M(u)A(u)}{u}y\right)^2=\frac{\alpha_\IF}{c^2(1-\alpha_\IF)^3}\frac{B}{2A}y^2=\frac{y^2}{2},
 \EQNY
 where $A$ and $B$ are given in (\ref{AB}).
Moreover, by (\ref{re})
$$\lim_{u\rw\IF}\frac{A(u)}{\Delta(u)}=\mathbb{Q}\lim_{u\rw\IF}\frac{\sigma(u)}{\inv{\sigma}(u^{-1}\sigma^2(u))}=\IF.$$
Hence \re{$y_2= \frac{y^2}{2} I_{\{y>0\}}-\frac{y^2}{2} I_{\{y<0\}}$}. Additionally, it follows from (\ref{y11})-(\ref{y12}) that $y_1=-\IF$. Thus  applying i) in  Theorem \ref{TH1}  we have
\BQN\label{pi6}
q_1(u)&\sim& \B_{X_{\varphi}}(x_2)\sqrt{\frac{2A\pi}{B}}\frac{1}{2\sqrt{\pi}}\int_{-\IF}^{y_2}|t|^{-1/2}e^{-|t|}dt\frac{u}{M(u)\Delta(u)}\Psi(M(u))\nonumber\\
&\sim&\B_{X_{\varphi}}(x_2)\sqrt{\frac{2A\pi}{B}}\Phi(y)\frac{u}{M(u)\Delta(u)}\Psi(M(u)), \quad u\rw\IF.
\EQN
Combination of (\ref{pi4upper}), (\ref{Bonfer4}) and (\ref{pi6}) leads to
\BQNY
q(u)\sim \B_{X_{\varphi}}(x_2)\sqrt{\frac{2A\pi}{B}}\Phi(y)\frac{u}{M(u)\Delta(u)}\Psi(M(u)), \quad u\rw\IF.
\EQNY
Consequently, by Theorem \ref{TH2} and (\ref{ratio})
\BQNY
\pk{\frac{\tau_u^*(x_1,x_2)-ut_u}{A(u)}<y}=\frac{\B_{X_{\varphi}}(x_2)}{\B_{X_{\varphi}}(x_1)}\Phi(y), \quad u\rw\IF
\EQNY
establishing the proof. \QED

\prooftheo{TH4} We set  $\Delta(u)=\overleftarrow{\sigma}\left(\frac{\sqrt{2}\sigma^2(T)}{u+cT}\right)$. Following (\ref{transfinite}), we have
\BQN\label{Bonfer5}
\lambda(u)\leq \pk{v(u)\int_0^T\mathbb{I}_u(X(t)-ct)\td t>x}\leq \lambda(u)+\lambda_1(u) ,
\EQN
where
\BQNY
\lambda(u)&=&\pk{v(u)\int_{E_{\re{4}}(u)}\mathbb{I}_0\left(\frac{X(t)}{u+ct}m(u)-m(u)\right)\td t>x},\\
\lambda_1(u)&=&\pk{\sup_{t\in [0,T]\setminus E_{\re{4}}(u)}\frac{X(t)}{u+ct}m(u)>m(u)},
\EQNY
with
$$E_{\re{4}}(u)=\left[T-\left(\frac{\ln m(u)}{ m(u)}\right)^2, T\right],
\quad v(u)=1/\Delta(u)$$ if
$\lim_{t\rw 0}\frac{|t|}{\sigma^2(|t|)}\in [0,\IF)$
and
 $$v(u)=(m(u))^2$$ if $\lim_{t\rw 0}\frac{|t|}{\sigma^2(|t|)}=\IF$. We shall derive the exact asymptotics of $\lambda(u)$ by applying Theorem \ref{TH1}, and show that
 $$\lambda_1(u)=o(\lambda(u)), \quad u\rw\IF.$$
We distinguish three cases: $\lim_{t\rw 0}\frac{|t|}{\sigma^2(|t|)}=0, (0,\IF)$ and $ \IF$, respectively.\\
{\it \underline{Case $\lim_{t\rw 0}\frac{|t|}{\sigma^2(|t|)}=0$ }}.\\
{\it \underline{Asymptotics of $\lambda(u)$}}.  Noting that $v(u)=1/\Delta(u)$ and scaling  time by $\Delta(u)$ we \re{get}
\BQN\label{pi70}
\lambda(u)=\pk{\int_{E_{\re{5}}(u)}\mathbb{I}_0\left(Z_u(t)-m(u)\right)\td t>x},
\EQN
with
\BQN\label{E5}
Z_u(t)=\frac{X(T-\Delta(u)t)}{u+cT-c\Delta(u)t}m(u), \quad E_{\re{5}}(u)=\left[0, (\Delta(u))^{-1}\left(\frac{\ln m(u)}{m(u)}\right)^2\right].
\EQN
In light of Lemma \ref{LemL2}, we have that
\BQN\label{varcor}
&&1-\sqrt{Var(Z_u(t))}\sim \frac{\dot{\sigma}(T)}{\sigma(T)}\Delta(u)|t|, t\in E_{\re{5}}(u), \nonumber\\
&&\lim_{u\rw\IF}\sup_{s\neq t, s,t\in E_{\re{5}}(u)}\left|\frac{\re{m^2(u)}(1-Corr(Z_u(t), Z_u(s)))}{\frac{\sigma^2(\Delta(u)|t-s|)}{\re{\sigma^2(\Delta(u))}}} -1\right|=0.
\EQN
With the notation \K{introduced in} Theorem \ref{TH1}, (\ref{varcor}) implies that
\BQN\label{notation2}
&&n(u)=m(u), \quad w=t, \quad g(u)=\frac{\dot{\sigma}(T)}{\sigma(T)}\Delta(u), \quad \Delta(u)=\overleftarrow{\sigma}\left(\frac{\sqrt{2}\sigma^2(T)}{u+cT}\right), \quad \eta=X,\nonumber\\
&&E(u)=[a_1(u), a_2(u)], \quad a_1(u)=0, \quad a_2(u)=(\Delta(u))^{-1}\left(\frac{\ln m(u)}{ m(u)}\right)^2.
\EQN
Next we check the conditions in i) of Theorem \ref{TH1}. It follows that
\BQN\label{check0}
\lim_{u\rw\IF}g(u)=\lim_{u\rw\IF}\frac{\dot{\sigma}(T)}{\sigma(T)}\Delta(u)=0, \quad g(u)a_1(u)=0, \quad \lim_{u\rw\IF}g(u)a_2(u)=\lim_{u\rw\IF}\left(\frac{\ln m(u)}{ m(u)}\right)^2=0,
\EQN
and
$$\gamma=\lim_{u\rw\IF}n^2(u)w(g(u))=\lim_{u\rw\IF}\frac{\dot{\sigma}(T)}{\sigma(T)}\Delta(u)(m(u))^2=
\mathbb{Q}\lim_{u\rw\IF}\frac{\overleftarrow{\sigma}(u^{-1})}{\sigma^2(\overleftarrow{\sigma}(u^{-1}))}=0.$$
Moreover,
$$x_1=\lim_{u\rw\IF}n^2(u)w(g(u)|a_1(u)|)=0,\quad x_2=\lim_{u\rw\IF}n^2(u)w(g(u)|a_2(u)|)=\re{\mathbb{Q}}\lim_{u\rw\IF}\frac{\dot{\sigma}(T)}{\sigma(T)}(\ln m(u))^2=\IF,$$
$$\lim_{u\rw\IF}a_2(u)=\re{\mathbb{Q}}\lim_{u\rw\IF}\frac{\sigma^2(\overleftarrow{\sigma}(u^{-1}))}{\overleftarrow{\sigma}(u^{-1})}\re{\ln m(u)}=\IF,$$
imply that $y_1=0$ and $y_2=\IF$. Additionally,
$$\lim_{u\rw\IF}n(u)(w(g(u)|a_1(u)|)+w(g(u)|a_2(u)|))=\lim_{u\rw\IF}\frac{\dot{\sigma}(T)}{\sigma(T)}\frac{(\ln m(u))^2}{m(u)}=0, \quad \lim_{u\rw\IF}\Delta(u)=0.$$
Consequently, by i) of Theorem \ref{TH1} and Remark \ref{remark} we have
\BQN\label{pi7}
\lambda(u)\sim \B_{B_{2\alpha_0}}(x)\frac{\sigma(T)}{\dot{\sigma}(T)}\frac{1}{(m(u))^2\Delta(u)}\Psi(m(u)).
\EQN
{\it \underline{Upper bound of $\lambda_1(u)$}}.  From \eqref{L3} \K{we have} \re{for sufficiently large $u$} that
$$\sup_{t\in [0,T]\setminus E_4(u)}Var\left(\frac{X(t)}{u+ct}m(u)\right)\leq 1-\mathbb{Q}\left(\frac{\ln( m(u))}{m(u)}\right)^2.$$
Further,  by {\bf BII}
\BQNY
(m(u))^2Var\left(\frac{X(t)}{u+ct}-\frac{X(s)}{u+cs}\right)&\leq& 2(m(u))^2\left(\frac{\sigma^2(|t-s|)}{(u+ct)^2}+\frac{\sigma^2(s)c^2(t-s)^2}{(u+ct)^2(u+cs)^2}\right)\\
&\leq& \mathbb{Q}\left(\sigma^2(|t-s|)+|t-s|^2\right)\leq \mathbb{Q}|t-s|^{\alpha_0}, \quad s,t \in [0,T].
\EQNY
Consequently, in light of  Piterbarg Theorem [\re{Theorem} 8.1 in \cite{Pit96}] we have that for $u$ sufficiently large
\BQN\label{pi8}
\lambda_1(u)\leq \mathbb{Q}_1 (m(u))^{2/\alpha_0} \Psi\left(\frac{m(u)}{\sqrt{1-\mathbb{Q}\left(\frac{\ln (m(u))}{m(u)}\right)^2}}\right).
\EQN
Combination of (\ref{Bonfer5}), (\ref{pi7}) and (\ref{pi8}) leads to
$$\pk{\frac{1}{\Delta(u)}\int_0^T\mathbb{I}_u(X(t)-ct)\td t>x}\sim  \B_{B_{2\alpha_0}}(x)\frac{\sigma(T)}{\dot{\sigma}(T)}\frac{1}{(m(u))^2\Delta(u)}\Psi(m(u)), $$
establishing  the claim.\\
{\it \underline{Case $\lim_{t\rw 0}\frac{|t|}{\sigma^2(|t|)}=1/\theta\in (0,\IF)$}}. First note that  (\ref{pi70})-(\ref{check0}) still hold. Next we check the conditions of ii) in Theorem \ref{TH1}. Following the notation in Theorem \ref{TH1}, we have that
\BQNY
&&\gamma=\lim_{u\rw\IF} n^2(u)w(g(u))=2\sigma(T)\dot{\sigma}(T)
\lim_{u\rw\IF}\frac{\overleftarrow{\sigma}(u^{-1})}{\sigma^2(\overleftarrow{\sigma}(u^{-1}))}
=\frac{2\sigma(T)\dot{\sigma}(T)}{\theta},\\
&&a_1=\lim_{u\rw\IF}a_1(u)=0, \quad a_2=\lim_{u\rw\IF} a_2(u)=\lim_{u\rw\IF} (\Delta(u))^{-1}\left(\frac{\ln m(u)}{ m(u)}\right)^2=\IF.
\EQNY
Consequently, in light of case ii) in Theorem \ref{TH1} and Remark \ref{remark}
\BQNY
\lambda(u)\sim \B_{B_{2\alpha_0}}^{\frac{2\sigma(T)\dot{\sigma}(T)|t|}{\theta}}(x)\Psi(m(u)), \quad u\rw\IF,
\EQNY
which together with (\ref{Bonfer5}) and (\ref{pi8}) establishes the claim.\\
{\it \underline{Case $\lim_{t\rw 0}\frac{|t|}{\sigma^2(|t|)}=\IF$ }}. Noting that $v(u)=(m(u))^2$,
scaling of time by $\Delta(u)$ we have
\BQNY
\lambda(u)=\pk{\int_{E_5(u)}\mathbb{I}_0\left(Z_u(t)-m(u)\right)\td t>\frac{x}{\Delta(u) (m(u))^2}},
\EQNY
with $ E_5(u)$ defined in (\ref{E5}) and
$$Z_u(t)=\frac{X(T-\Delta(u)t)}{u+cT-c\Delta(u)t}m(u).$$
Next we verify conditions of case iii) in Theorem \ref{TH1}.
Notice that  (\ref{pi70})-(\ref{check0}) still hold. Using notation in Theorem \ref{TH1}, we have
$$\gamma=\lim_{u\rw\IF} n^2(u)w(g(u))=2\sigma(T)\dot{\sigma}(T)
\lim_{u\rw\IF}\frac{\overleftarrow{\sigma}(u^{-1})}{\sigma^2(\overleftarrow{\sigma}(u^{-1}))}=\IF,$$  $$b_1=\lim_{u\rw\IF}\frac{a_1(u)g(u)}{\inv{w}(n^{-2}(u))}=0,
b_2=\lim_{u\rw\IF}\frac{b(u)g(u)}{\inv{w}(n^{-2}(u))}=\lim_{u\rw\IF}\frac{\dot{\sigma}(T)}{\sigma(T)}\frac{(\ln m(u))^2}{m(u)}=\IF.$$
Moreover,$$ \frac{x}{\Delta(u) (m(u))^2}=\frac{\inv{w}(n^{-2}(u))}{g(u)}\frac{\dot{\sigma}(T)}{\sigma(T)}x =\theta(u)\frac{\dot{\sigma}(T)}{\sigma(T)}x.
$$
Consequently, by case iii) in Theorem \ref{TH1} and Remark \ref{remark}
\BQNY
\lambda(u)=\pk{\int_{E_5(u)}\mathbb{I}_0\left(Z_u(t)-m(u)\right)\td t>\frac{\inv{w}(n^{-2}(u))}{g(u)}\frac{\dot{\sigma}(T)}{\sigma(T)}x}\sim \mathcal{B}_0^{|t|}\left(\frac{\dot{\sigma}(T)}{\sigma(T)}x, [0,\IF)\right)\Psi(m(u)), \quad u\rw\IF,
\EQNY
which combined with (\ref{Bonfer5}), (\ref{pi8}) and (\ref{Piterbarg}) yields that
$$\pk{(m(u))^2\int_0^T\mathbb{I}_u(X(t)-ct)\td t>x}\sim e^{-\frac{\dot{\sigma}(T)}{\sigma(T)}x}\Psi(m(u)), \quad u\rw\IF. $$
This completes the proof. \QED
\proofkorr{TH5}. Observe that for any $y\geq 0$
\BQN
\pk{\frac{\dot{\sigma}(T)}{\sigma^3(T)}u^2(T-\tau_{u,T}^*(x_1,x_2))\geq y}=\frac{\pk{\tau_{u,T}(x_2)\leq T_u(y)}}{\pk{\tau_{u,T}(x_1)\leq T}},
\EQN
with $T_u(y)=T-u^{-2}\frac{\sigma^3(T)}{\dot{\sigma}(T)}y$. By the definition of $\tau_{u,T}(x)$, we have
\BQN\label{taux2}
\pk{\tau_{u,T}(x_1)\leq T}&=&\pk{v(u)\int_0^T\mathbb{I}_u(X(t)-ct)\td t>x_1},\nonumber\\
 \pk{\tau_{u,T}(x_2)\leq T_u(y)}&=&\pk{v(u)\int_0^{T_u(y)}\mathbb{I}_u(X(t)-ct)\td t>x_2},
\EQN
with $v(u)=1/\overleftarrow{\sigma}\left(\frac{\sqrt{2}\sigma^2(T)}{u+cT}\right)$ if $\lim_{t\rw 0}\frac{|t|}{\sigma^2(|t|)}\in [0,\IF)$, and $v(u)=(m(u))^2$ if $\lim_{t\rw 0}\frac{|t|}{\sigma^2(|t|)}=\IF$.
Since $\limit{u}T_u(y)=T$, we get that the \re{asymptotics} of (\ref{taux2}) is the same as in Theorem \ref{TH4} with $T$ replaced by $T_{u}(y)$ and $x$ replaced by $x_2$. Using this fact and  by Theorem \ref{TH4} for all $\lim_{t\rw 0}\frac{|t|}{\sigma^2(|t|)}\in [0,\IF]$
\BQNY
\frac{\pk{\tau_{u,T}(x_2)\leq T_u(y)}}{\pk{\tau_{u,T}(x_1)\leq T}}&\sim& \Gamma(x_1,x_2)\frac{\Psi\left(\frac{u+cT_u(y)}{\sigma(T_u(y))}\right)}{\Psi(m(u))}\\
&\sim& \Gamma(x_1,x_2)\exp\left( -\frac{1}{2}\left(\frac{u+cT}{\sigma(T)}\right)^2\left(\left(1-\frac{c(T-T_u(y)) }{u+cT}\right)^2\left(1-\frac{\dot{\sigma}(T)(T-T_u(y))}{\sigma(T)}\right)^{-2}-1\right)\right)\\
&\sim& \Gamma(x_1,x_2)e^{-y}, \quad y\geq 0,
\EQNY
where $\Gamma(x_1,x_2)$ is given in Theorem \ref{TH5}.
This completes the proof.\QED\\

\prooftheo{selfsimilar0} By self-similarity of $X$, we have that
\BQNY
\pk{\frac{1}{u\overleftarrow{\rho}((\widehat{A}u^{1-H})^{-2})}\int_0^\IF \mathbb{I}_u(X(t)-ct)dt>x}&=&\pk{\frac{1}{\overleftarrow{\rho}((\widehat{A}u^{1-H})^{-2})}\int_0^\IF \mathbb{I}_0\left(\frac{X(t)}{1+ct}-u^{1-H}\right)dt>x}\\
&=&\pk{\frac{1}{\overleftarrow{\rho}((\widehat{A}u^{1-H})^{-2})}\int_0^\IF \mathbb{I}_0\left(\widehat{A}\frac{X(t)}{1+ct}-\widehat{A}u^{1-H}\right)dt>x}.
\EQNY
Let  \BQN\label{zt}
Z(t)=\widehat{A}\frac{X(t)}{1+ct},\quad n(u)=\widehat{A}u^{1-H}, \quad  \Delta(u)=\overleftarrow{\rho}((n(u))^{-2}).
\EQN
Observe that
\BQN\label{piselfsimilar}
\varpi(u)\leq \pk{\frac{1}{\Delta(u)}\int_0^\IF \mathbb{I}_0(Z(t)-n(u))dt>x}\leq \varpi(u)+\varpi_1(u),
\EQN
where
\BQNY
\varpi(u)&=&\pk{\frac{1}{\Delta(u)}\int_{[t_0-(\ln n(u))/n(u), t_0+(\ln n(u))/n(u)]} \mathbb{I}_u(Z(t)-n(u))dt>x},\\
\varpi_1(u)&=& \pk{\sup_{t\in [0,\IF)\setminus [t_0-(\ln n(u))/n(u), t_0+(\ln n(u))/n(u)]}Z(t)>n(u)}.
\EQNY
First we derive the asymptotics of $\varpi(u)$ by applying Theorem \ref{TH1} and then show that $\omega_1$ is asymptotically negligible in comparison to $\varpi(u)$ as $u\rw\IF$.
We distinguish between three cases: $\gamma=0, \gamma\in (0,\IF)$ and $\gamma=\IF$. \\
{\it \underline{Case $\gamma=0$}}.\\
{\it \underline{The asymptotics of $\varpi(u)$}}. In order to apply  Theorem \ref{TH1}, we rewrite
\BQN\label{pi9}
\varpi(u)&=&\pk{\int_{[-\delta_u, \delta_u]} \mathbb{I}_0(Z(t_0+\Delta(u)t)-n(u))dt>x},
\EQN
with $\delta_u=\frac{\ln n(u)}{n(u)\Delta(u)}$. By (\ref{selfcor}) and (\ref{selfvariance}), we have that, as $u\rw\IF$
\BQN\label{var11}
1-\sqrt{Var(Z(t_0+\Delta(u)t))}\sim \frac{\widehat{B}}{2\widehat{A}}(\Delta(u)t)^2, \quad  t\in [-\delta_u, \delta_u],
\EQN
and
\BQNY
\lim_{u\rw\IF}\sup_{s\neq t, s,t \in [-\delta_u, \delta_u]}\left|\frac{n^2(u)(1-Corr(Z(t_0+\Delta(u)s), Z(t_0+\Delta(u)t)))}{\frac{\rho(\Delta(u)|t-s|)}{\rho(\Delta(u))}}-1\right|=0.
\EQNY
Thus we have that, corresponding to the notation in Theorem  \ref{TH1},
\BQN\label{notationth2.1}
w(t)=\frac{\widehat{B}}{2\widehat{A}}t^2, \quad g(u)=\Delta(u), \quad a_1(u)=-\delta_u, \quad a_2(u)=\delta_u.
\EQN
A direct calculation shows that
\BQN\label{gamma}
\lim_{u\rw\IF}n^2(u)w(g(u))=\frac{\widehat{B}}{2\widehat{A}}\lim_{u\rw\IF}(n(u)\Delta(u))^2 =\frac{\widehat{B}}{2\widehat{A}}\lim_{u\rw\IF}(n(u)\overleftarrow{\rho}((n(u))^{-2}))^2=\frac{\widehat{B}}{2\widehat{A}}\lim_{t\rw 0}\frac{t^2}{\rho(|t|)}=0.
\EQN
\BQNY
\lim_{u\rw\IF}n(u)w(g(u)|a_i(u)|)=\frac{\widehat{B}}{2\widehat{A}}\lim_{u\rw\IF}\frac{(\ln n(u))^2}{n(u)}=0.
\EQNY
\BQNY
x_i=\lim_{u\rw\IF} n^2(u)w(g(u)|a_i(u)|)=\frac{\widehat{B}}{2\widehat{A}}\lim_{u\rw\IF}(\ln n(u))^2=\IF, \quad i=1,2.
\EQNY
By the fact that
$$\lim_{u\rw\IF} a_1(u)=-\lim_{u\rw\IF} \frac{\ln n(u)}{n(u)\Delta(u)}=-\IF, \quad \lim_{u\rw\IF}a_2(u)=\IF,$$
we have that
$$y_1=-\IF, \quad y_2=\IF.$$
Noting that $\lim_{u\rw\IF} \Delta(u)=0$ and by i) in  Theorem \ref{TH1} and Remark \ref{remark} i), we have
\BQNY
\varpi(u)&\sim& \mathcal{B}_{B_\alpha}(x)2^{-1}\int_{-\IF}^\IF |t|^{-1/2}e^{-|t|}dt\frac{\sqrt{2\widehat{A}/\widehat{B}}}{n(u)\Delta(u)}\Psi(n(u))\\
&\sim& \B_{B_\alpha}(x)\sqrt{\frac{2\widehat{A}\pi}{\widehat{B}}} \frac{1}{\overleftarrow{\rho}((\widehat{A}u^{1-H})^{-2})\widehat{A}u^{1-H}}\Psi(\widehat{A}u^{1-H}).
\EQNY
{\it \underline{Upper bound of $\varpi_1(u)$}}. Observe that, for $u$ sufficiently large,
\BQNY
\varpi_1(u)\leq \sum_{i=2}^{4}\varpi_{i}(u),
\EQNY
where
$$\varpi_2(u)=\pk{\sup_{t\in [0,\epsilon]}Z(t)>n(u)}, \quad \varpi_3(u)=\pk{\sup_{t\in [\epsilon, L]\setminus [t_0-(\ln n(u))/n(u), t_0+(\ln n(u))/n(u)]}Z(t)>n(u)},$$
 $$\varpi_4(u)=\sum_{k=0}^\IF\pk{\sup_{t\in [L+k, L+k+1]}Z(t)>n(u)},$$
 with $L>t_0$.
 In order to prove that $\varpi_1(u)=o(\varpi(u)), u\rw\IF$,  it suffices to show that $\varpi_i(u), i=2, 3, 4$ are negligible compared with $\varpi(u)$ respectively. We begin with $\varpi_2(u)$. Using the fact that for $\epsilon>0$ sufficiently small there exists $0<\delta<1$ such that
 $$\sup_{t\in [0,\epsilon]} Var(Z(t))<1-\delta,$$
 and by Borell-TIS inequality (see \cite{GennaBorell, AdlerTaylor}), we have
 \BQNY
 \varpi_2(u)\leq e^{-\frac{\left(n(u)-\E{\sup_{t\in [0,\epsilon]}Z(t)}\right)^2}{2(1-\delta)}}=o(\varpi(u)), \quad u\rw\IF.
 \EQNY
 Next we focus on $\varpi_3(u)$. By (\ref{selfcor}) and self-similarity of $X$, we have that for $s,t\geq \epsilon>0$ and $r=\frac{t_0}{t}$
 \BQNY
 \E{\left(\overline{X}(t)-\overline{X}(s)\right)^2}&=&2(1-Corr(X(s), X(t)))\\
 &=&2(1-Corr(X(rs), X(rt)))\\
 &\leq& 4\rho(r|s-t|)\leq 8\left(t_0/\epsilon\right)^{\alpha/2}|t-s|^{\alpha/2}, \quad |s-t|\rw 0,
 \EQNY
 which indicates that for $s,t\geq \epsilon>0$ and $|s-t|\leq L_1<\IF$ with $L_1$ a positive constant
 \BQN\label{holder}
  \E{\left(\overline{X}(t)-\overline{X}(s)\right)^2}\leq \mathbb{Q}|t-s|^{\alpha/2}.
 \EQN
 Hence (\ref{selfvariance}) and Piterbarg inequality (Lemma 8.1 in \cite{Pit96}) leads to
 \BQNY
 \varpi_3(u)&\leq& \pk{\sup_{t\in [\epsilon, 2t_0]}\overline{X}(t)>\frac{n(u)}{1-\mathbb{Q}_1(\ln n(u))^2/n^2(u)}}\\
 &\leq& \mathbb{Q} (n(u))^{4/\alpha} \Psi\left(\frac{n(u)}{1-\mathbb{Q}_1(\ln n(u))^2/n^2(u)}\right)\\
 &=& o(\varpi(u)), \quad u\rw\IF.
 \EQNY
 Finally, we consider $\varpi_4(u)$. Using the fact that
 $$\sup_{t\in [L+k, L+k+1]} Var\left(Z(t)\right)\leq \left(\frac{\widehat{A}}{c}(L+k)^{H-1}\right)^2$$
 and  (\ref{holder}), by Piterbarg inequality we have that, for $L$ sufficiently large,
 \BQNY
\varpi_4(u)&\leq& \sum_{k=0}^\IF\pk{\sup_{t\in [L+k, L+k+1]}\overline{X}(t)>\frac{n(u)}{\frac{\widehat{A}}{c}(L+k)^{H-1}}}\\
 &\leq& \sum_{k=0}^\IF\mathbb{Q}(n(u))^{4/\alpha} \Psi\left(\frac{n(u)}{\frac{\widehat{A}}{c}(L+k)^{H-1}}\right)\\
 &\leq & \mathbb{Q}(n(u))^{4/\alpha}\Psi\left(\frac{n(u)}{\mathbb{Q}_1L^{H-1}}\right)=o(\varpi(u)), \quad u\rw\IF.
 \EQNY
 Consequently, we conclude that
 \BQN\label{negpi10}
\varpi_1(u)=o(\varpi(u)), \quad u\rw\IF.
 \EQN
 Therefore,
 \BQNY
 \pk{\frac{1}{u\overleftarrow{\rho}((\widehat{A}u^{1-H})^{-2})}\int_0^\IF \mathbb{I}_u(X(t)-ct)\td t>x}\sim \B_{B_\alpha}(x)\sqrt{\frac{2\widehat{A}\pi}{\widehat{B}}} \frac{1}{\overleftarrow{\rho}((\widehat{A}u^{1-H})^{-2})\widehat{A}u^{1-H}}\Psi(\widehat{A}u^{1-H}).
 \EQNY
 This establishes the claim for $\gamma=0$.\\
 {\it \underline{Case $\gamma\in (0,\IF)$}}.\\
 {\it \underline{The asymptotics of $\varpi(u)$}}. Note that (\ref{pi9})-(\ref{notationth2.1}) also hold for $\gamma\in (0,\IF)$. As in (\ref{gamma}), we have that
 \BQNY
 \lim_{u\rw\IF}n^2(u)w(g(u))=\frac{\widehat{B}}{2\widehat{A}}\lim_{t\rw 0}\frac{t^2}{\rho(|t|)}=\frac{\widehat{B}}{2\widehat{A}}\gamma\in (0,\IF).
 \EQNY
 Moreover,
 $$\lim_{u\rw\IF}g(u)a_i(u)=0, \quad\lim_{u\rw\IF} a_1(u)=-\lim_{u\rw\IF} \frac{\ln n(u)}{n(u)\Delta(u)}=-\IF, \quad \lim_{u\rw\IF}a_2(u)=\IF.$$
 Hence by ii) in Theorem \ref{TH1}, we have that
 \BQNY
 \varpi(u)\sim\widehat{\B}_{B_2}^{\frac{\widehat{B}\gamma t^2}{2\widehat{A}}}(x)\Psi(\widehat{A}u^{1-H}),
 \EQNY
which combined with  (\ref{negpi10}) and (\ref{Pit2}) establishes the claim.\\
{\it \underline{Case $\gamma=\IF$}}.\\
  {\it \underline{The asymptotics of $\varpi(u)$}}. Using the notation for $Z(t), \Delta(u), n(u),\delta_u, \varpi_1(u)$ in (\ref{zt})-(\ref{pi9}), we have that
   \BQNY
\varpi_5(u)\leq \pk{\int_0^\IF \mathbb{I}_u(X(t)-ct)dt>xu^H}\leq \varpi_1(u)+\varpi_5(u),
   \EQNY
   where
   $$\varpi_5(u)=\pk{\int_{[-\delta_u, \delta_u]} \mathbb{I}_0(Z(t_0+\Delta(u)t)-n(u))dt>\frac{xu^{H-1}}{\Delta(u)}}.$$
  We focus on $\varpi_5(u)$. Note that in this case, (\ref{var11})-(\ref{notationth2.1}) still hold. Following the notation in Theorem \ref{TH1}, we have that
    \BQNY
 \lim_{u\rw\IF}n^2(u)w(g(u))=\frac{\widehat{B}}{2\widehat{A}}\lim_{t\rw 0}\frac{t^2}{\rho(|t|)}=\IF, \quad \lim_{u\rw\IF} g(u)|a_i(u)|=\lim_{u\rw\IF}\frac{\ln n(u)}{n(u)}=0, i=1,2.
 \EQNY
 Moreover,
 $$b_1=\lim_{u\rw\IF}\frac{a_1(u)g(u)}{{\inv{w}(n^{-2}(u))}}=-\sqrt{\frac{\widehat{B}}{2\widehat{A}}}\lim_{u\rw\IF}\ln n(u)=-\IF, \quad b_2=\lim_{u\rw\IF}\frac{a_2(u)g(u)}{{\inv{w}(n^{-2}(u))}}=\sqrt{\frac{\widehat{B}}{2\widehat{A}}}\lim_{u\rw\IF}\ln n(u)=\IF.$$
   In this case, using the notation in (\ref{notationth2.1}) we have $$\frac{xu^{H-1}}{\Delta(u)}=x\sqrt{\frac{\widehat{A}\widehat{B}}{2}}\inv{w}(n^{-2}(u))/g(u).$$

 In light of iii) in Theorem \ref{TH1}, we have that
 \BQNY
 \varpi_5(u)&=&\pk{\int_{[-\delta_u, \delta_u]} \mathbb{I}_0(Z(t_0+\Delta(u)t)-n(u))dt>x\sqrt{\frac{\widehat{A}\widehat{B}}{2}}\inv{w}(n^{-2}(u))/g(u)}\\
 &\sim& \mathcal{B}_0^{t^2}\left(\sqrt{\frac{\widehat{A}\widehat{B}}{2}}x, (-\IF, \IF)\right)\Psi(n(u)),
 \EQNY
 which together with (\ref{negpi10}) and (\ref{Piterbarg}) completes the proof.
  \QED

 \prooftheo{selfsimilar1} It follows that
 \BQNY
 \pk{v(u)\int_0^T \mathbb{I}_u(X(t)-ct)\td t>x}=\pk{v(u)\int_0^T \mathbb{I}_0\left(\frac{X(t)}{u+ct}\frac{u+cT}{T^{2H}}-\frac{u+cT}{T^{2H}}\right)\td t>x},
 \EQNY
 where
 $$X_u(t)=\frac{X(t)}{u+ct}\frac{u+cT}{T^{2H}}.$$
Direct calculation shows that
\BQNY
\lim_{u\rw\IF}\lim_{t\rw T, t<T}\left|\frac{1-Var\left(X_u(t)\right)}{|T-t|}-\frac{H}{T}\right|=0.
\EQNY
Moreover, by {\bf S} with $t_0=T$,
\BQNY
\lim_{\epsilon\rw 0}\sup_{s\neq t, T-\epsilon\leq  s, t\leq T}\left|\frac{1-Corr\left(X_u(t), X_u(s)\right)}{\rho(|t-s|)}-1\right|=0.
\EQNY
The local behavior of variance and correlation functions is the same as the one in Theorem \ref{TH4} (see Lemma \ref{LemL2}). Hence
 proceeding similarly as in the proof of Theorem \ref{TH4}, that is replacing $\frac{\dot{\sigma}(T)}{\sigma(T)}$ by $\frac{H}{T}$,  $\frac{\sigma^2(\cdot)}{2\sigma^2(T)}$ by $\rho(\cdot)$ and $\sigma(T)$ by $T^H$ in the proof of Theorem \ref{TH4}, we establish the claim.\QED

\section{Appendix}

{\bf Proof of (\ref{upper2})}. In order to prove (\ref{upper2}), it suffices to prove that for $c_1>0$
$$\int_{0}^{g(u)a_2(u)}e^{-c_1n^2(u)w(t)}\td t\sim c_1^{-1/\beta} \inv{w}(n^{-2}(u))\beta^{-1}\int_{0}^{c_1y_2} t^{1/\beta-1}e^{-t}\td t, \quad u\rw\IF.$$
Recall that $w$ is a regularly varying function at $0$ with index $\beta>0$ satisfying
$$\lim_{u\rw\IF}n^2(u)w(g(u)a_2(u))=y_2>0,$$
where  $a_2(u)>0$ and $\lim_{u\rw\IF} g(u)=\IF$ and $\lim_{u\rw\IF}g(u)a_2(u)=0$.
We can assume that $w(x)= \ell(x) x^\beta$ with $\ell$ normalized slowly varying function at 0. Then from \cite{BI1989}, we know that
$\ell(x) x^\beta$ is ultimately monotone for any $\beta\not=0$, $\ell$ is continuously differentiable and
\BQN\label{ldif}
 \lim_{x\to0}  \frac{x \ell'(x)}{\ell(x)}=0.
 \EQN
Let $c_1n^2(u)w(t)=s$. Then we have that
\BQNY
\int_{0}^{g(u)a_2(u)}e^{-c_1n^2(u)w(t)}\td t=\frac{1}{c_1n^2(u)}\int_{0}^{c_1n^2(u)w(g(u)a_2(u))}
\frac{1}{w'\left(\overleftarrow{w}\left(\frac{s}{c_1n^2(u)}\right)\right)}e^{-s}\td s.
\EQNY
By (\ref{ldif}), it follows that
$$w'(x)\sim \beta \frac{w(x)}{x}, \quad x\rw 0.$$
Hence
\BQNY
\int_{0}^{g(u)a_2(u)}e^{c_1n^2(u)w(t)}\td t&\sim& \int_{0}^{c_1n^2(u)w(g(u)a_2(u))}
\frac{\overleftarrow{w}\left(\frac{s}{c_1n^2(u)}\right)}{\beta s}e^{-s}\td s\\
&=& \frac{\overleftarrow{w}\left(\frac{1}{c_1n^2(u)}\right)}{\beta}\int_{0}^{c_1n^2(u)w(g(u)a_2(u))}
\frac{\overleftarrow{w}\left(\frac{s}{c_1n^2(u)}\right)}{\overleftarrow{w}\left(\frac{1}{c_1n^2(u)}\right)}s^{-1}e^{-s}\td s.
\EQNY
Note that $\overleftarrow{w}$ is a regularly varying function at $0$ with index $1/\beta$.  Moreover, using Potter's bound for $\overleftarrow{w}$  (\re{see e.g.}
\cite{BI1989} \re{or} \cite{Soulier} or Lemma 6.1 in \cite{KEP20151}), we have that for any $\epsilon\in (0, \min(1,1/\beta))$ and all $u$ large
\BQNY
\frac{\overleftarrow{w}\left(\frac{s}{c_1n^2(u)}\right)}{\overleftarrow{w}\left(\frac{1}{c_1n^2(u)}\right)}\leq
(1+\epsilon) s^{1/\beta-\epsilon}, \quad 0<s<c_1n^2(u)w(g(u)a_2(u)).
\EQNY
Moreover, {since $\overleftarrow{w}$ is regularly varying at $0$}, then for any $s>0$
$$\lim_{u\rw\IF}\frac{\overleftarrow{w}\left(\frac{s}{c_1n^2(u)}\right)}{\overleftarrow{w}
\left(\frac{1}{c_1n^2(u)}\right)}=s^{1/\beta}.$$
Hence {the dominated convergence theorem implies} that
$$\int_{0}^{g(u)a_2(u)}e^{-c_1n^2(u)w(t)}\td t\sim c_1^{-1/\beta} \inv{w}(n^{-2}(u))\beta^{-1}\int_{0}^{c_1y_2} s^{1/\beta-1}e^{-s}\td s.$$
This completes the proof.
\QED

{\bf Acknowledgement}:
The authors would like to thank Enkelejd Hashorva
for his comments and suggestions that significantly improved
the content of this contribution.
The authors also would like to thank the referees for their valuable comments and remarks.
K. D\c ebicki
was partially supported by NCN Grant No 2015/17/B/ST1/01102 (2016-2019).
Financial support from the Swiss National Science
Foundation Grant 200021-175752 is also kindly acknowledged.

\bibliographystyle{ieeetr}

\bibliography{queue2d}

\end{document}